\documentclass[12 pt]{article}
\usepackage{amssymb,amsmath,amsfonts,amsthm,xcolor}
\newcommand\BBP{{\mathbb {P}}}

\newcommand\bkR{{\mathbb {R}}}
\newcommand\BBN{{\mathbb {N}}}

\newcommand\bkE{{\mathbb {E}}}

\newcommand\N{{\mathbb {N}}}
\newcommand\R{{\mathbb {R}}}

\newtheorem {Theorem}{Theorem}[section]
\newtheorem {Lemma}[Theorem]{Lemma}
\newtheorem {Proposition}[Theorem]{Proposition}

\theoremstyle{definition}
\newtheorem{Definition}{Definition}[section]
\newtheorem{Notation}{Notation}[section]
\newtheorem{Remark}{Remark}[section]

\newtheorem{nota}{Notation}[section]

\renewcommand{\P}{\ensuremath{\mathbb {P}}}
\newcommand{\E}{\mathbb{E}}

\newcommand{\F}{\mathcal{F}}
\newcommand{\G}{{\mathcal G}}

\newcommand{\as}{\mbox{$\P$-a.s.}}

\newcommand\beq{\begin{equation}}
\newcommand\eeq{\end{equation}}

\def\ti{\tilde}

\def\as{\mbox{$\P$-a.s.}}
\def\mfa{{\mathfrak a}}

\def\B{{\mathcal B}}
\def\f{{\mathfrak f}}
\def\mfa{{\mathfrak a}}

\def\B{{\mathcal B}}
\def\f{{\mathfrak f}}
\def\g{{\mathfrak g}}
\def\k{{\mathfrak k}}
\def\s{{\mathfrak s}}

\def\PP{{\mathcal P}}
\def\u{{\mathfrak u}}

\begin{document}

\title{Rates of convergence  in invariance principles for  random walks on linear groups via martingale methods}

\author{
    C.~Cuny\footnote{Universit\'e de Brest, LMBA, UMR CNRS 6205.
    Email:~christophe.cuny@univ-brest.fr},
    J.~Dedecker\footnote{Universit\'e Paris Descartes, Sorbonne Paris Cit\'e,  Laboratoire MAP5 (UMR 8145).
        \newline
    Email:~jerome.dedecker@parisdescartes.fr},
    Florence Merlev\`ede\footnote{Universit\'{e} Paris-Est, LAMA (UMR 8050), UPEM, CNRS, UPEC.
    Email:~florence.merlevede@u-pem.fr}
}

\maketitle

\begin{abstract} In this paper, we give explicit rates in the central limit theorem and  in the almost sure invariance principle for general ${\mathbb R}^d$-valued cocycles that appear in the study of the  left  random walk on linear groups. Our method of proof lies on a suitable martingale approximation and on a careful estimation of some coupling coefficients linked with the underlying Markov structure. Concerning the martingale part, the available results in the literature are not accurate enough to give almost optimal rates whether in the central limit theorem for  the Wasserstein distance, or in  the strong approximation.  A part of this paper is devoted to circumvent this issue. We then exhibit  near optimal rates both in the central limit theorem in terms of Wasserstein distance and in the  almost sure invariance principle for  ${\mathbb R}^d$-valued martingales with stationary increments having moments of order $p \in ]2, 3]$ (the case of sequences of reversed martingale differences is also considered). Note also that, as an application of our results for general ${\mathbb R}^d$-valued cocycles,  a special attention is paid to the Iwasawa cocycle and the Cartan projection for reductive Lie groups. 
\end{abstract}


\smallskip


\section{Introduction}
\setcounter{equation}{0}

Let $G=GL_d(\R)$ be the group of invertible real matrices of order $d\ge 2$. Denote by $\|\cdot\|$ the operator norm on $G$ associated with the euclidean norm on $\R^d$.  Let $\mu$ be a probability on the Borel sets of $G$. Let $(Y_n)_{n\ge 1}$ be independent identically distributed (iid) variables with law $\mu$. It is well-known that, if  $\mu$ admits a moment of order 1 and  under some irreducibility assumption on the support of $\mu$, then the sequence
$(\frac1n \log \|Y_n\cdots Y_1\|)_{n\ge 1}$ converges almost surely  to $\lambda_\mu\in \R$   (see for instance Furstenberg and Kesten \cite{FK}).

 \medskip

If $\mu$ is further assumed to admit  exponential moments and to satisfy a  proximality property, 
Le Page \cite{L} and   Guivarc'h and Raugi \cite{GR} proved a central limit theorem (CLT) with rate as well as other probabilistic results, by  spectral gap methods. The CLT has been obtained by Jan 
\cite{Jan} under a moment of order $2+\varepsilon$, by mean of martingale approximation. 

\medskip

Very recently, Benoist and Quint \cite{BQ} managed to prove the CLT under a moment of order 2 thanks to an ingenious and explicit martingale-coboundary 
decomposition. Later, Cuny, Dedecker and Jan \cite{CDJ} gave a different proof of Benoist-Quint's result. Using precise controls of some coupling coefficients of the underlying Markov chain and an explicit martingale approximation, they derived  more probabilistic results such as almost sure invariance principles with rates. Next, Cuny, Dedecker and 
Merlev\`ede \cite{CDM} obtained different types of deviation results. 

\medskip

In all the above mentioned works, the first step of the proof consists in reducing the study to a suitable cocycle on $G \times X$, where $X$ stands for the
projective space of $\R^d$. Recall that $\sigma\, :G\times X\to \R$ is a cocycle if $\sigma(gg',x)=\sigma(g,g'\cdot x)+\sigma(g',x)$ for every $g,g'\in G$ and $x\in X$, where $\cdot$ denotes an action of $G$ on $X$. For instance to deal with $(\log  \|Y_n\cdots Y_1\|)_{n \geq 1}$, one can use the cocycle $ \displaystyle \sigma(g,x) = \log  \Big (  \frac{ \| g \cdot x\|}{ \| x\|}\Big ) $. 

\medskip

The central limit theorem for cocycles benefited from an active research in the last years. Let us mention, beside \cite{BQ}, \cite{CDJ} and \cite{CDM}, the works of Bj\"orklund \cite{bjorklund}, Benoist and Quint \cite{BQ2} and Horbez \cite{Horbez}. 

\medskip

It happens that the arguments developped in \cite{CDJ} are somewhat general and apply equally to any cocycle on $\sigma \, : G\times X\to \R^d$ ($G$ a locally compact group acting on  a compact metric space  $X$) provided that the action is suitably contracting and that the cocycle is Lipschitz in the second coordinate with Lipschitz norm $\sigma_{{\rm Lip}}(g)$ satisfying integrability conditions with respect to $\mu$, see Definition 3.3 and  Section 3.2 for more details.

\medskip

Hence a first goal of the paper is to extend results  of \cite{CDJ} to general $\R^d$-valued  cocycles. As a main motivation we have in mind the Iwasawa cocycle to which our results apply and for which the CLT has been obtained by Benoist and Quint \cite{BQ} under a moment of order 2. Then, as in \cite{CDJ} or \cite{CDM} many probabilistic results follow: strong laws of large numbers with rate, the CLT (and its functional form), deviation inequalities. For those results it is not an issue to treat vector-valued processes.  However, there are several results in probability theory that do not extend easily from the 
one-dimensional case to the multivariate one. As a matter of fact, for $\R^d$-valued martingales, no almost sure invariance principle (ASIP) with \emph{explicit} rate 
 is available in the litterature. Indeed, the only known rates are of the form $O(n^{1/2-\varepsilon})$ for some 
$\varepsilon>0$, see for  instance Morrow and Philipp \cite{MorP}, 
Monrad and Philipp \cite{MonP} and Eberlein \cite{Eberlein}, while in the one dimensional case,  martingales with stationary increments  in $L^p$ for some  for $2<p<4$, satisfy an ASIP with  rate $o(n^{1/p}\sqrt{\log n})$. A similar rate also holds for $p=4$ but no better rate is available for $p>4$.

\medskip

A second goal of the paper is then to prove an ASIP with explicit rates for $\R^d$-valued martingales with stationary  increments, which is the main technical result of the present work. Our result hold 
when the martingales are in $L^p$ for some $2<p\le 3$, so that the case $3<p\le 4$ remains open, contrary to the one-dimensional case. 
The method of proof also allows to derive  rates of convergence in the CLT for $\R^d$-valued martingales with stationary increment, in terms of the Wasserstein distance 
$W_1$.  Using then a suitable martingale approximation, we derive ASIP with explicit rates for general $\R^d$-valued cocycles. 

\medskip

Let us mention that, to prove the ASIP for cocycles, a second way would be to apply a multidimensional version of the strong approximation result of Berkes, Liu and Wu \cite{BLW}, as given in Karmakar and Wu \cite{KW}. 
However, a direct application of this result would not give the good rate of convergence with respect to the moment of $\mu$ (see the introduction of 
\cite{CDM2} for more details). Hence an adaptation of the result of \cite{KW}, similar to what we did in \cite{CDM2} for real-valued cocycles, would be necessary here, and it is not clear 
wether this adaptation is  feasible or not. Let us also mention that  the ASIP for ${\mathbb R}^d$-valued martingales is interesting  in itself, and can be useful in other situations.

\medskip

The paper is organized as follows. In Section 2 we state our results for $\R^d$-valued martingales. In section 3, we describe the type of \emph{weakly contractive} actions that we shall deal with, and obtain several preliminary results. Then, we describe the properties that the cocycles should satisfy to implement the arguments used in \cite{CDJ}, and  we derive several probabilistic results for those cocycles, emphasizing the consequences of  our martingale 
results. We conclude Section 3 by an  application  to the Iwasawa cocycle and to the Cartan projection. The fact that  the Iwasawa cocycle does satisfy our conditions follows from the recent book by Benoist and Quint \cite{BQ-book}. Finally, Section 4 is devoted to the proofs and the appendix contains a useful Fuk-Nagaev type inequality for martingales.

\bigskip

\section{Rates of convergence  in the ASIP and the CLT  for ${\mathbb R}^d$-valued martingales}\label{Sec:mart}

\setcounter{equation}{0}

In this section we consider  a stationary sequence of random variables with values in ${\mathbb R}^d$ and defined on a probability
space $(\Omega,{\mathcal{A}},{\mathbb{P}})$. Let $(\F_n)_{n \in {\mathbb N}}$ be a stationary and non-decreasing sequence of $\sigma$-algebras (see e.g. page 10 in \cite{MPU} for the definition of a stationary filtration). We suppose that $(d_n)_{n \in {\mathbb N}}$ is a stationary sequence of  martingale differences with respect to $(\F_n)_{n \in {\mathbb N}}$, i.e. for each $n$, $d_n$ is integrable, $\F_n$-measurable and such that  $\E(d_n|\F_{n-1})=0$ $\P$-a.s. 

For any random variable $X$ with values in ${\mathbb R}^d$, we shall use the notation $(X)_i$ to mean its $i$-th coordinate. In addition $|\cdot |_d$ means the euclidean norm on ${\mathbb R}^d$.

In this section we give, under
projective conditions, rates of
convergence in the almost sure invariance principle for the partial sums associated with a stationary sequence of  martingale differences with values in ${\mathbb R}^d$.

\begin{Notation}\label{notquant}
For any $p > 2$, define the envelope norm $\Vert \, . \, \Vert_{1,
\Phi, p}$ by
$$
\Vert X \Vert_{1, \Phi, p}  = \int_0^1 (1 \vee \Phi^{-1} (1-u/2)
)^{p-2} Q_X(u) du
$$
where  $\Phi$ denotes the d.f. of the $N(0, 1)$ law, and $Q_X$
denotes the quantile function of $|X|$, that is the cadlag inverse
of the tail function $x \rightarrow {\mathbb P}(|X|>x)$. 
%Note that for any $\gamma >0$, there exists a positive constant 
%$C(p, \gamma)$ such that $\Vert X \Vert_{1, \Phi, p} \leq C(p, \gamma) \Vert X \Vert_{1 + \gamma}$. 
\end{Notation}
\begin{Remark} \label{comparison-norms}Let $a>1$ and $p>2$.
Applying H\"older's inequality, we see that there exists a positive
constant $C(p, a)$ such that $\Vert X \Vert_{1, \Phi, p} \leq C(p,a)
\|X\|_a$.
\end{Remark}

\begin{Theorem} \label{Thmart} Let $(d_n)_{n \in {\mathbb Z}}$ be a ${\mathbb R}^d$-valued stationary sequence of  martingale differences   
with respect to $(\F_n)_{n \in {\mathbb Z}}$. Let $M_n = \sum_{k=1}^n d_k$. Let $p \in ]2,3]$. Assume that $\bkE |d_0 |_d^p < \infty$ and that for any $i,j$ such that $1 \leq i,j \leq d$,
 \beq \label{C1}
\sum_{n=1}^{\infty}\frac{1}{n^{3-p/2}}\big \Vert \bkE \big
((M_n)_i (M_n)_j\big | \F_{0} \big ) - \bkE \big
((M_n)_i (M_n)_j \big ) \big \Vert_{1, \Phi, p} < \infty \, ,\eeq and \beq \label{C2}
\sum_{n=1}^{\infty}\frac{1}{n^{1+2/p}}\big \Vert \bkE \big
((M_n)_i (M_n)_j\big | \F_{0} \big ) - \bkE \big
((M_n)_i (M_n)_j \big ) \big \Vert_{p/2} < \infty \, .\eeq Then, 
\begin{enumerate}
\item For any $\varepsilon >0$, enlarging $\Omega$ if
necessary, there exists a sequence $(N_i)_{i
\geq 1}$ of iid ${\mathbb R}^d$-valued centered gaussian random variables with ${\rm Var} (N_1) = \E (d_0 d_0^t)$ such that 
\[
 M_n - \sum_{i=1}^n N_i  =  \left\{
  \begin{aligned}
  o ( n^{1/p} (\log n)^{ \frac{p+1}{2p} + \varepsilon})& \text{ if $p\in ]2,3[ $}  \\
 o ( n^{1/3} (\log n)^{1 + \varepsilon}) & \text{ if $p=3 $} \\
  \end{aligned}
\right. \text{ almost surely} \, .
\]
\item Enlarging $\Omega$ if
necessary, there exists a sequence $(N _i)_{i
\geq 1}$ of iid ${\mathbb R}^d$-valued centered gaussian random variables with ${\rm Var} (N_1) = \E (d_0 d_0^t)$ such that 
\[
 \Big \Vert  \sup_{1 \leq k \leq n} \big |   M_k - \sum_{i=1}^k N_i  \big |_d   \Big \Vert_1=  \left\{
  \begin{aligned}
  O ( n^{1/p} (\log n)^{ \frac{p-1}{2p} })& \text{ if $p\in ]2,3[ $}  \\
 O ( n^{1/3} (\log n)^{2/3}) & \text{ if $p=3 $} \\
  \end{aligned}
\right. \, .
\]
\end{enumerate}
\end{Theorem}
\begin{Remark}\label{twoconds} According to Remark \ref{comparison-norms}, if $p \in ]2, 3]$,  condition (\ref{C1}) is implied by the condition: for some $a>1$,
 \beq \label{C1alter}
\sum_{n=1}^{\infty}\frac{1}{n^{3-p/2}}\big \Vert \bkE \big
((M_n)_i (M_n)_j\big | \F_{0} \big ) - \bkE \big
((M_n)_i (M_n)_j \big ) \big \Vert_{a} < \infty\, .\eeq
Hence, both \eqref{C1} and \eqref{C2} hold as soon as  \beq 
\sum_{n=1}^{\infty}\frac{1}{n^{3-p/2}}\big \Vert \bkE \big
((M_n)_i (M_n)_j\big | \F_{0} \big ) - \bkE \big
((M_n)_i (M_n)_j \big ) \big \Vert_{{p/2}} < \infty\, .\eeq
\end{Remark}

Let  ${\rm Lip}( |\cdot|_d)$ be the set of Lipschitz functions $g$ from ${\mathbb R}^d$ to ${\mathbb R}$ such that 
$|g(x) - g(y)| \leq |x-y|_d$. For two measures $\mu$ and $\nu$ on ${\mathbb R}^d$, let 
\[
W_1 ( \mu, \nu) = \sup_{g \in {\rm Lip}( |\cdot|_d)}  (    \mu (g) - \nu (g)  )  \, .
\]

\begin{Theorem} \label{ThmartW1cond} Let $(d_n)_{n \in {\mathbb Z}}$ be a ${\mathbb R}^d$-valued stationary sequence of  martingale differences   
with respect to $(\F_n)_{n \in {\mathbb Z}}$ satisfying the assumptions of Theorem \ref{Thmart}. Let $G_{\Sigma}$ be a ${\mathbb R}^d$-valued centered gaussian random variables with ${\rm Var} (G_\Sigma) =  \E (d_0 d_0^t):= \Sigma$. Then 
\[
\left \Vert  W_{1} \left (P_{n^{-1/2} M_n | {\mathcal F}_0}, P_{G_{\Sigma}}  \right ) \right  \Vert_1 \leq  \left\{
  \begin{aligned}
  &  C  n^{ (2-p) /2} & \text{ if $p \in ]2,3[$} \\
 &  C n^{-1/2}  \log n  & \text{ if $p=3 $}   \\
  \end{aligned}
\right. \, ,
\]
where $C$ is a positive constant depending on $(p,d)$ but not on $n$.
\end{Theorem}

\begin{Remark} \label{remark-ThmartW1cond} 
It follows from Theorem \ref{ThmartW1cond} that, for any bounded ${\mathcal F}_0$-measurable random variable 
$Z$, 
$$
    \sup_{g \in {\rm Lip}( |\cdot|_d)} \left ( {\mathbb E} \left ( Z g \left  ( n^{-1/2} M_n \right ) \right ) - {\mathbb E}(Z) {\mathbb E}(g(G_\Sigma)) \right )
     \leq  \left\{
    \begin{aligned}
  &  C  n^{ (2-p) /2} & \text{ if $p \in ]2,3[$} \\
 &  C n^{-1/2}  \log n  & \text{ if $p=3 $} \\
  \end{aligned}
\right. \, .
$$
In particular, taking $Z \equiv 1$, we obtain that 
\[
 W_{1} \left (P_{n^{-1/2} M_n}, P_{G_{\Sigma}}  \right )  \leq  \left\{
  \begin{aligned}
  &  C  n^{ (2-p) /2} & \text{ if $p \in ]2,3[$} \\
 &  C n^{-1/2}  \log n  & \text{ if $p=3 $} \\
  \end{aligned}
\right. \, .
\]
From the above result, proceeding as in \cite{DMR09}, one can derive rates in the central limit theorem with respect to $W_1$ for the normalized partial sums of a large class of stationary sequence of ${\mathbb R}^d$-valued dependent random variables satisfying some mixingale type conditions. For related results under other types of dependence  conditions, let us mention the paper by P\`ene \cite{Pe}. 
\end{Remark}
\begin{Remark}[Reversed martingale differences sequences] \label{commentreversed} Let $p \in ]2,3]$.  Assume that $(d_n)_{n \in {\mathbb Z}}$ is a ${\mathbb R}^d$-valued stationary sequence of {\it   reversed } martingale differences  in ${\mathbb L}^p$ with respect to a stationary and non-increasing sequence $({\mathcal G}_n)_{n \in {\mathbb Z}}$ of $\sigma$-algebras. This means that for any integer $n$, $d_n$ is ${\mathcal G}_n$-adapted and ${\mathbb E} ( d_n | {\mathcal G}_{n+1} ) = 0$ a.s. Let $M_n = \sum_{k=1}^n d_k$. We infer that 
the conclusions of Theorem \ref{Thmart} and Theorem \ref{ThmartW1cond}  hold for $M_n$ (with   ${\mathcal G}_{n+1} $ in place of  $ {\mathcal F}_0$ in Theorem \ref{ThmartW1cond}) provided that the conditions \eqref{C1} and \eqref{C2}  are replaced by the following ones: 
 \beq \label{C1rev}
\sum_{n=1}^{\infty}\frac{1}{n^{3-p/2}}\big \Vert \bkE \big
((M_n)_i (M_n)_j\big | \G_{n+1} \big ) - \bkE \big
((M_n)_i (M_n)_j \big ) \big \Vert_{1, \Phi, p} < \infty \, ,\eeq and \beq \label{C2rev}
\sum_{n=1}^{\infty}\frac{1}{n^{1+2/p}}\big \Vert \bkE \big
((M_n)_i (M_n)_j\big | \G_{n+1} \big ) - \bkE \big
((M_n)_i (M_n)_j \big ) \big \Vert_{p/2} < \infty \, .\eeq 
See Section \ref{proofcomment} for a short proof of this remark.

Proceeding, as in \cite{CM15}, this type of result for reversed martingale differences sequences allows to derive rates in the strong invariance principle for ${\mathbb R^d}$-valued observables of a large class of dynamical systems (see also \cite{KKM} where reversed martingale approximations in ${\mathbb R^d}$ are provided that could be used to verify conditions such as \eqref{C1rev} and \eqref{C2rev}). This reversed martingale approximation method allows to derive rates up to $n^{1/3}$ (up to additional logarithm terms) which improves, in case of  ${\mathbb R^d}$-valued bounded H\"older observables, the rates obtained in \cite{MN09} which depend on the dimension $d$ (in particular, when $d$ is large, the rates in \cite{MN09} are close to $1/2$). Note that, for a class of processes with exponential decay of correlations, the rate $n^{1/4 + \varepsilon}$ for some $\varepsilon >0$ can be reached for ${\mathbb R^d}$-valued observables, by  using the method developed in \cite{Go10}. 
\end{Remark}

\section{Weakly contracting actions}\label{section-cocycle}
\setcounter{equation}{0}

\subsection{Definitions and properties}
We shall work in the general situation described in section 11 of \cite{BQ-book}. 
At first reading one may assume that the group $F$ below is trivial (i.e. reduced to its 
neutral element). In our applications, $F$ will be trivial as soon as $G$ is connected (see Remark \ref{connected}).   

\smallskip

Let $G$ be a locally compact second countable group. Let $s\, :\, 
G \to F$ be a continuous homomorphism onto a finite group $F$.

\smallskip 

Let $\mu$ be a probability on the Borel sets $\B(G)$ of $G$. Let 
$s(\mu):=\mu\circ s^{-1}$,  and $\Delta_\mu$ the support of $s(\mu)$. 

\begin{Definition}
We say that $\mu$  is \emph{$F$-adapted} if $\Delta_\mu$ spans 
$F$ as a subgroup. We say that $\mu$ is {$F$-strictly aperiodic} if it is $F$-adapted and if the smallest normal subgroup of $F$ a class of which contains $\Delta_\mu$ is $F$.
\end{Definition}
\begin{Remark}
Our terminology is different from \cite{BQ-book}. Our definition of 
{$F$-strict aperiodicity} corresponds to their definition of 
\emph{aperiodicity} (see Derriennic-Lin \cite[Prop. 1.6]{DL} for a proof of the equivalence of the definitions). 
%The 
%reason why we use a different terminology is that 
%in the theory of random walks, several  authors use \emph{aperiodic} 
%to mean \emph{adapted}. 
\end{Remark}
\begin{Remark}\label{connected}
When $F$ is trivial, $\mu$ is automatically $F$-adapted and $F$-strictly aperiodic. If $G$ is connected then $s(G)=\{e_F\}$ ($e_F$ is the neutral element of $F$). 
In particular, in that case, if $\mu$ is $F$-adapted, $F$ must be trivial.
\end{Remark} 

Let $X$ be a  compact 
and second countable metric space on which $G$ acts continuously (we denote that action by $g\cdot x$). Assume 
that $X$ is fibered over $F$, meaning that 
there exists a $G$-equivariant continuous mapping $f\, :\, X\to F$. 
Recall that $G$-equivariance means that $f(g\cdot x)=s(g)f( x)$ for every 
$(g,x)\in G\times X$.

\smallskip Notice that when $X$ is fibered over $F$ then the sets $(X_\f)_{\f\in F}:=(f^{-1}(\{\f\}))_{\f\in F}$ are compact  and open and obviously form a partition of $F$. Those sets are called the fibers.

\medskip

Before going further, to motivate our definitions, let us describe some situations to which our results will apply. Take $G= GL_d(\R)$. Then, it is well known 
 that we have an homeomorphism $G= KN$ where $K=O_d(\R)$ is the group of orthogonal matrices of size $d$ and $N$ is the group of upper triangular matrices of size $d$  with positive entries on the diagonal. We take $X:= G/N$ which is homeomorphic to $K$ (hence compact). Then, if $x=k N\in 
X$, we set $f(x)= {\rm sgn\, det} ( k)$, the sign of the determinant of $k$. Similarly, for $g\in G$ we set $s(g)={\rm sgn\, det} ( g)$. The fact that the conditions imposed below may  be satisfied will be explained in Section \ref{Iwasawa-section}.

One may also take $G:= SL_d(\R)$. Then, $K=SO_d(\R)$ and $f$ and $s$ become trivial.

When $G=GL_d(\R)$ or $SL_d(\R)$, one may also take for $X$ the projective space of $\R^d$. This special case has been handled in \cite{CDJ} and \cite{CDM}.

%Let $\mu$ be a probability measure on the Borel sets of $G$.  

\begin{Definition} 
We say that a probability $\nu$
on the Borel sets $\B(X)$ of $X$ is $\mu$-invariant if
$$
\int_{G\times X}h(g\cdot u)\mu (dg) \nu(du)=
\int_{ X}h(u) \nu(du)\, ,
$$
for every bounded Borel function $h$ on $X$.
\end{Definition}

For every $g\in G$, set

\begin{equation}\label{lipschitz}
{\rm Lip}(g):=\sup_{ x\neq y, f(x)=f(y)} \frac{d(g\cdot x,g\cdot y)}{d(x,y)}
\, ,
\end{equation} 
where the supremum is over all $x,y\in X$ such that $x\neq y$ and $f(x)=f(y)$.

\medskip

\medskip

\begin{Definition}\label{weakly-contracting}
Let $p\ge 1$. We say that the action of $G$ on $(X,d)$ is \emph{weakly} $(\mu,p)$-\emph{contracting}
if 
\begin{equation}\label{moment-lip}
\int_G \left(\log^+ ({\rm Lip}(g))\right)^p\, \mu(dg)<\infty\, ,
\end{equation}
where $\log^+ (x) = \log (\max (1, x) )$, 
and if there exists $n_0\in \N$, such that 
\begin{equation}\label{contracting}
\sup_{ x\neq y,f(x)=f(y)} \int_G \log \Big( \frac{d(g\cdot x,g\cdot y)}{d(x,y)}
\Big)\, \mu^{*n_0}(dg) <0\, .
\end{equation}
When $p=1$ we shall just say that the action is weakly $\mu$-contracting.
\end{Definition}
\begin{Remark}\label{comparewithBQ}
Notice that the left-hand side in \eqref{contracting} may be $-\infty$.
Benoist and Quint \cite[Definition 11.1]{BQ-book} called an  action 
$(\mu,\gamma)$-\emph{contracting} for some $\gamma>0$, if 
\begin{equation}\label{moment-lip-BQ}
\int_G  ({\rm Lip}(g))^\gamma \, \mu(dg)<\infty\, ,
\end{equation}
and 
\begin{equation}\label{contracting-BQ}
\sup_{ x\neq y,f(x)=f(y)} \int_G \Big( \frac{d(g\cdot x,g\cdot y)}{d(x,y)}
\Big)^\gamma\, \mu^{*n_0}(dg) <1\, .
\end{equation}
Using Jensen's inequality, one can see that \eqref{moment-lip-BQ} implies 
\eqref{moment-lip} and that \eqref{contracting-BQ} implies \eqref{contracting}.
\end{Remark}
\begin{Remark}\label{rem-holder}
Let $\alpha \in (0,1]$, then $d_\alpha(x,y):=(d(x,y))^\alpha$ defines another 
metric on $X$. Let $p\ge 1$. For every $\alpha\in (0,1]$, the action of $G$ on 
$(X, d)$ is weakly $(\mu,p)$-contracting, for some  $p\ge 1$ if only if 
the action of $G$ on $(X,d_\alpha)$ is. 
\end{Remark}

\medskip

Benoist and Quint \cite{BQ-book} proved (see their Lemma 11.5 p. 171) 
that if the action is  $(\mu,\gamma)$-contracting, for some $\gamma>0$,
then there exists a unique $\mu$-invariant probability on $\B(X)$. 

\smallskip

We shall prove that their result remains true under the weaker assumption that 
the action is weakly $\mu$-contracting. As in \cite{BQ-book}, the proof makes use 
of the left (and right) random walk on $X$.

\smallskip

\medskip

Let $(\Omega,\F,\P)$ be a probability space. Assume that there exists 
$(Y_n)_{n\ge 1}$ iid random variables on $(\Omega,\F,\P)$ taking values 
in $G$ with common law $\mu$. Define $A_n:=Y_n\cdots Y_1$ and $B_n= 
Y_1\cdots Y_n$ for every      
$n\ge 1$ and $A_0={\bf e}$, where ${\bf e}$ stands for the neutral element of $G$.

\medskip

\begin{Proposition}\label{unique}
Assume that the action is weakly $\mu$-contracting and that 
$\mu$ is $F$-adapted. Then, there exists a unique 
$\mu$-invariant probability on $\B(X)$. 
\end{Proposition}
%\noindent {\bf Proof.} 

The proof is based on the following two lemmas. The proofs of the proposition and of the lemmas below are postponed to Section \ref{sectproofsection-cocycle}.

\begin{Lemma}\label{basic-lemma}
Assume that $\mu$ is $F$-adapted. Let $\nu$ be a $\mu$-invariant probability 
on $\B(X)$. Then, $\nu(X_f)=\frac1{|F|}$ for every $f\in F$.
\end{Lemma}

\begin{Lemma}\label{lemma-complete}
Assume that the action is weakly $\mu$-contracting. Then, there exists $\ell >0$ such that  \begin{equation}\label{complete-int}
\sum_{n\ge 0} \sup_{x\neq y, f(x)=f(y)}\P\, \left( \log\, \Big(\frac{d(B_{2^n}\cdot x,B_{2^n}\cdot 
y)}{d(x,y)}\Big)\ge -2^n\ell\right)\, <\infty\, .
\end{equation}
\end{Lemma}
\begin{Remark}
Of course, since $(A_n)_{n\ge 1}$ has the same law as $(B_n)_{n\ge 1}$, the lemma holds with $A_{2^n}$ in place of $B_{2^n}$. However, we shall need the above form in the proof of Proposition \ref{unique}.
\end{Remark}

\smallskip

We shall see now that for weakly-contracting actions, the trajectories (starting from a same fiber) of the 
left random walk on $X$ are almost surely arbitrary close. 
The next lemma is a version of Lemma 6 in  \cite{CDJ}.

\begin{Lemma} \label{lem-comp}
Let $p\ge 1$. Assume that the action is weakly $(\mu,p)$-contracting. 
Then, there exists $\ell>0$, such that 
\begin{equation}\label{complete-cocycle}
\sum_{k\ge 1} k^{p-2}\max_{k\le j\le 2k}\sup_{ x\neq y,f(x)=f(y)} \P \left(
 \log \left(d(A_{j-1}\cdot  x, 
A_{j-1}\cdot  y)\right)\ge -\ell k\right)<\infty\, .
\end{equation}
Moreover, there exists $\delta>0$ such that for every $x,y\in X$, with $f(x)=f(y)$, 
 \begin{equation}\label{exp-conv}
 d(A_n\cdot x,A_n\cdot y)= O({\rm e}^{-\delta n})\qquad 
 \as
\end{equation}
\end{Lemma}

\subsection{Cocycles over weakly contracting actions}\label{cocycle}

Our goal is to obtain limit theorems  for \emph{cocyles} associated with a weakly 
$\mu$-contracting $G$-action on a compact metric space  $X$ that is fibered over 
$F$, as in the previous section. We shall be concerned with cocycles taking values in a finite dimensional $\R$-vector  space $E$, equipped with an euclidean norm $|\cdot |_E$. 

\begin{Definition}
We say that $\sigma\,:\, G\times X\to E$ is a cocycle if 
for every $g,g'\in G$ and every $u\in X$,
\begin{equation}\label{cocycle-property}
\sigma(gg',u)=\sigma(g,g'\cdot u) + \sigma(g',u)\, .
\end{equation}
\end{Definition}
%\begin{Remark}
%Given an orthonormal basis $(e_i^*)_{1\le i\le m}$ of $E^*$, $e_i^*\circ 
%\sigma$ is a real-valued cocycle for every $1\le i\le m$, hence in many of 
%the results below it is enough to consider the case where $E=\R$. 
%\end{Remark}

\smallskip

Of course we need some regularity assumptions on the considered cocycles.  
In order to state the needed assumptions we shall introduce some notations.

\smallskip

Given a cocycle $\sigma \, :\, G\times X\to E$, define, for every $g\in G$,
\begin{gather}
\sigma_{\rm sup}(g):= \sup_{x\in X}|\sigma(g,x) |_E\, ,\\
\sigma_{\rm Lip}(g):= \sup_{x\neq y, f(x)=f(y)} 
\frac{|\sigma(g,x)-\sigma(g,y)|_E}{d(x,y)}\\
{\rm and }\qquad \kappa_0(g):= \max\big(\sigma_{\rm sup}(g),\log(\sigma_{\rm Lip}(g)
\big) \, .
\end{gather}

\begin{Definition}
We say that a cocycle $\sigma\, : \, G\times X\to E$ has a polynomial 
moment of order $p\ge 1$ if 
\begin{equation}
 \int_G \kappa_0^p(g)\, \mu(dg)
<\infty \, .
\end{equation}
\end{Definition}
%\begin{Remark}
%Assume that the action of $G$ on $(X,d)$ is weakly  $(\mu,p)$-contracting of order $p\ge 1$. 
%In the sequel, we shall study cocycles having polynomial 
%moments of order $p$ (with the same $p\ge 1$). In particular, we need to control the Lipschitz constant $\sigma_{\rm Lip}$. Using Remark \ref{rem-holder}, one can 
%see that using the distance $d_\alpha$ rather than $d$ it is enough to control the 
%H\"older constant of order $\alpha$: 
%\[
%\sigma_\alpha (g):= \sup_{x\neq y, f(x)=f(y)} 
%\frac{|\sigma(g,x)-\sigma(g,y)|_E}{(d(x,y))^\alpha} \, .
%\]
%\end{Remark}
\medskip

\begin{Remark}\label{remark-cocycle}
Let $\sigma \, :\, G\times X\to E$ be a cocycle. Let $(e_i)_{1\le i\le d}$ be an orthonormal basis of $E$. Then, for every $1\le i\le d$, 
$(g,x)\mapsto \sigma_i(g,x):= \langle e_i, \sigma(g,x)\rangle$ is an $\R$-valued cocycle and it is not hard to see that $\sigma$ has a polynomial moment of order $p\ge 1$ if and only if $\sigma_i$ has a polynomial moment of order $p\ge 1$ for every $1\le i\le d$. Because of this it will be sufficient to deal with $\R$-valued cocycles. 
\end{Remark}
\textit{Hence, due to the remark above,  we assume that the dimension $d=1$ in the remaining of this subsection.}

\medskip

We first recall the strong law of large numbers for cocyles. This  follows here from Theorem 3.9 of \cite{BQ-book}, using Proposition \ref{unique}. 

\begin{Proposition}\label{SLLN}
Assume that the action is weakly $\mu$-contracting and that $\mu$ is 
$F$-adapted. Assume moreover that $\sigma_{\rm sup}\in L^1(\mu)$.
 Then, for every $x\in X$,
$$
\frac{\sigma(A_n,x)}{n}\underset{n\to \infty}\longrightarrow 
\int_{G\times X} \sigma(g,x)\mu(dg) \nu(dx) =:\lambda_\mu \, ,
$$
$\P$-a.s. and in $L^1(\P)$. Moreover, the convergence in $L^1(\P)$ is uniform 
over $x\in X$.
\end{Proposition}
%\noindent {\bf Proof.} Using Proposition  \ref{unique}, the result 
%\hfill $\square$

\medskip 

For every $q>0$, define a non decreasing, concave function 
$H_q$ on $[0,1]$ by $H_q(0)=0$ and for every $x\in (0,1]$,
$H_q(x)= \big|\log (x{\rm e}^{-q-1})\big|^{-q}$.

\medskip

The next result is a version of Lemma 5 of \cite{CDJ}. The proof being identical, it  is therefore omitted.
%\marginpar{la preuve etant identique a CDJ, peut etre peut-on l'omettre}

\begin{Lemma}\label{lem-reg}
For every $r>1$, there exists 
$C_r>0$ such that for every $g\in G$ and every $x,y \in X$, 
with $f(x)=f(y)$,
\begin{equation}\label{regularity}
|\sigma(g, x)-\sigma(g,  y)|\le C_r (1+
\kappa_0(g))^r \,  H_{r-1}(d(x, y))\, .
\end{equation}
\end{Lemma}
%\noindent   {\bf Proof.} By definition  for every $g\in G$ and  every $ x, y\in X$, 
%with $f(x)=f(y)$
% \begin{equation}\label{BQ-est}
% |\sigma(g, x)- \sigma(g, y)|\le {\rm e}^{\kappa_0(g)}d(x, y)\, ,
% \end{equation}
% and
%\begin{equation}\label{trivial-est}
%|\sigma(g,x)|\le \kappa_0(g)) \, .
%\end{equation}

% Assume that $d( x, y)\le {\rm e}^{-\kappa_0(g)}$. Using that 
%$t \mapsto t(H_{r-1}(t))^{-1}$ 
%is non decreasing on  $(0, {\rm e}]$, we have 
%$${\rm e}^{\kappa_0(g)} d( x, y)\le   
%\frac{H_{r-1}(d( x, y))}{H_{r-1}({\rm e}^{-\kappa_0(g)})}\, .
%$$
%Hence, by \eqref{BQ-est}, 
%\begin{equation}\label{first-est}
%|\sigma(g,- \sigma(g, y)|\le \frac{H_{r-1}(d( x, y))}
%{H_{r-1}({\rm e}^{-\kappa_0(g)}} 
%\le (r+\kappa_0(g) )^{r-1}H_{r-1}(d( x, y))\, .
%\end{equation}

%Assume now that $ d( x, y)>{\rm e}^{-\kappa_0(g)}$. By \eqref{trivial-est}, 

%\begin{equation}\label{second-est}
%|\sigma(g, x)- \sigma(g, y)|\le 2\frac{\kappa_0(g)
% H_{r-1}({\rm e}^{-\kappa_0(g) })}
%{H_{r-1}({\rm e}^{-\kappa_0(g) })}\le 2 (r+\kappa_0(g))^{r}
%H_{r-1}(d( x, y))
%\,.
%\end{equation}

%Combining \eqref{first-est} and \eqref{second-est}, we see that 
%\eqref{regularity} holds. \hfill $\square$

Lemma \ref{lem-reg} allows us to prove that the almost sure behaviour 
of the cocycle along the random walk does not depend on the starting point of a given fiber. 

\begin{Lemma}\label{startingpoint-0}
 Let $p>1$. Assume that the action is weakly $(\mu,p)$-contracting, that 
$\mu$ is $F$-adapted and that 
$\sigma$ admits a moment of order $p$. Then, for every $x,y\in X$ 
such that $f(x)=f(y)$,
\begin{equation*}
|\sigma(A_n,x)-\sigma(A_n,y)|=O(1)\qquad \as
\end{equation*}
\end{Lemma}
\noindent {\bf Proof.} Since $\sigma$ admits a moment of order $p$, 
it is a well-known consequence of the Borel-Cantelli lemma that 
$\kappa_0 (Y_k)=o(k^{1/p})$ \as\,  Let $x,y\in X$ be such that 
$f(x)=f(y)$. Let $r>2p/(p-1)$. Using the cocycle property,  \eqref{regularity} and \eqref{exp-conv}, we have, $\P$-a.s.
\begin{gather*}
|\sigma(A_n,x)-\sigma(A_n,y)| \le \sum_{k=1}^n |\sigma(Y_k,A_{k-1}\cdot 
x)- \sigma(Y_k,A_{k-1}\cdot y)|\\ \le C_r(1+\kappa_0(Y_k))^r 
\sum_{k=1}^n H_{r-1}(d(A_{k-1}\cdot x,A_{k-1}\cdot y) ) 
\le \tilde C_r\sum_{k=1}^n  \frac{k^{r/p}}{k^{r-1}}\leq \sum_{k\ge1 }\frac{\tilde C_r}{k^{r(p-1)/p-1}}<\infty
\, , 
\end{gather*} 
and the result follows.

\begin{Proposition}\label{cond-cobord}
Let $p>1$. Assume that the action is weakly $(\mu,p)$-contracting, that 
$\mu$ is $F$-adapted and that 
$\sigma$ admits a moment of order $p$. Then, for $q\in[1,p)$,
\begin{equation}\label{decadix}
\sum_{k=1}^\infty k^{p-q-1} 
\sup_{ f(x)=f(y)} 
\E \left (|\sigma(Y_k,A_{k-1}\cdot x)-\sigma(Y_k,A_{k-1}\cdot y)|^{q} 
\right)< \infty \, \, ,
\end{equation}
and for $q \in (0,1]$, 
\begin{equation}\label{decadix2}
\sum_{k=1}^\infty k^{p-2} 
\sup_{f(x)=f(y)} 
\E \left (|\sigma(Y_k,A_{k-1}\cdot x)-\sigma(Y_k,A_{k-1}\cdot y)|^{q} \right)  < \infty \, \, .
\end{equation}
\end{Proposition}
The proof of the proposition may be done as the proof of Proposition 3 of 
\cite{CDJ}, using Lemmas \ref{lem-comp} et \ref{lem-reg}.

Proposition \ref{cond-cobord} gives us a control along trajectories starting from a same fiber. To deal with trajectories starting from arbitrary fibers we shall need the following lemma (see for instance \cite[Lemma 11.6b]{BQ-book}). 

\begin{Lemma}\label{expo-lemma}
Assume that $\mu$ is $F$-strictly aperiodic. There exist $C>0$ and $0<\rho<1$ such that for every bounded function $\varphi$ on $F$, 
$$
 \sup_{\f \in F}\left|\E\big( \varphi(s(A_n)^{-1}\f\big)-\frac1{|F|}\sum_{\f'\in F} 
 \varphi(\f')\right|\le C\rho^n \sup_{\f \in F}|\varphi(\f)|
 \qquad \forall n\in \N\, .
$$
\end{Lemma}
\begin{Proposition}\label{prop-aperiodic}
Let $p>1$. Assume that the action is weakly $(\mu,p)$-contracting, that 
$\mu$ is $F$-strictly aperiodic and that 
$\sigma$ admits a moment of order $p$. Then, 
\begin{equation}\label{simple-series}
 \sum_{k\ge 1} k^{p-2}\sup_{ x \in X}|\E(\sigma(Y_k,A_{k-1}
 \cdot x))-\lambda_\mu| <\infty \, .
 \end{equation}
\end{Proposition}
\noindent {\bf Proof.} One easily sees that $\big(\sup_{ x \in X}\big|\E(\sigma(Y_k,A_{k-1}\cdot x))- \lambda_\mu \big|\big)_{k\ge 1}$ is non increasing. Hence, it 
is enough to prove that 
\begin{equation*}
 \sum_{k\ge 1} k^{p-2}\sup_{ x \in X}|\E(\sigma(Y_{2k},A_{2k-1}
 \cdot x))-\lambda_\mu| <\infty \, .
 \end{equation*}

We shall make use of the following identity based on Lemma \ref{basic-lemma} 
$$
\sigma(Y_{2k},A_{2k-1}\cdot x)= |F|\sum_{\f \in F} 
\int_{X_\f} \nu(dy)\sigma(Y_{2k},A_{2k-1}\cdot x)\, {\bf 1}_{
X_\f}(A_{k-1}\cdot x)\, .
$$
Notice that $A_{k-1}\cdot x\in X_\f$ if and only if $s(A_{k-1})^{-1}\f
=f(x)$. Using Lemma \ref{expo-lemma} and independence, we infer that
\begin{align*}
&  \left|\E(\sigma(Y_{2k},A_{2k-1}\cdot x))-\lambda_\mu \right| \\ 
&  \quad  =
\left| \E\left(|F|\sum_{\f \in F} 
\int_{X_\f}\Big( \sigma(Y_{2k},A_{2k-1}\cdot x)- \sigma(Y_{2k},Y_{2k-1}\ldots Y_k\cdot y)\Big)\, {\bf 1}_{
X_\f}(A_{k-1}\cdot x)\nu(dy)\right)\right.\\
& \quad  \quad \quad \left.\qquad \qquad + |F|\sum_{\f \in F} \int_{ X_\f} \E
\big(\sigma(Y_{2k},Y_{2k-1}\ldots Y_k\cdot y)\big)
\nu(dy)\, \left(\E \left( {\bf 1}_{
\{f(x)\}}(s(A_{k-1})^{-1}\f ) \right)- \frac1{|F|}\right)\right|\\
&  \quad  \le \sup_{ y\, :\, f(y)=f(x)} 
\E \left (|\sigma(Y_k,A_{k-1}\cdot x)-\sigma(Y_k,A_{k-1}\cdot y)|
\right)
 +C \rho^{k-1}|F|\int_G\sigma_{\rm sup}(g)\mu(dg)\, .
\end{align*}
Then, the result follows from \eqref{decadix} with $q=1$. \hfill $\square$ 
\medskip

It is also possible to prove the following version of Proposition 4 of 
\cite{CDJ}. Define for every integer $j$ and every $x\in X$, $\tilde X_j(x):= \sigma(Y_j,A_{j-1}\cdot x)-\lambda_\mu$.

\begin{Proposition}\label{newprop2}
 Let $p>2$. Assume that the action is weakly $(\mu,p)$-contracting, that 
 $\mu$ is $F$-adapted  and that 
 $\sigma$ admits a moment of order $p$. Then
 \begin{equation}
 \label{square-series}\sum_{k\ge 1} k^{p-3}
 \sup_{  f(x)=f(y)}
 \E \left (\left|\tilde X^2_n(x) -\tilde X^2_n(y) \right | \right)<\infty\, ,
 \end{equation}
 and for every $\gamma<p-3+1/p$,
 \begin{equation}
 \label{mixed-series}\sum_{k\ge 1} k^{\gamma}\sup_{ f(x)=f(  y)}\sup_{ k
 \le j<i \le 2k}\E \left (
 \left |\tilde X_i(x)\tilde X_j(x)-
 \tilde X_i(y)\tilde X_j(y)
 \right | \right)<\infty\, .
 \end{equation}
 \end{Proposition}
\begin{Proposition}\label{newprop2-aperiodic}
Let $p>2$. Assume that the action is weakly $(\mu,p)$-contracting, that 
$\mu$ is $F$-strictly aperiodic and that 
$\sigma$ admits a moment of order $p$. Then,
 \begin{equation}
 \label{square-series2}\sum_{k\ge 1} k^{p-3}
 \sup_{ x  \in X}\left| \E\, (
 \tilde X_k^2(x)) -\int_{X}\E(\tilde X_k^2(y))\nu(dy)\right | <\infty\, ,
 \end{equation}
 and for every $\gamma<p-3+1/p$,
 \begin{equation}
 \label{mixed-series2}\sum_{k\ge 1} k^{\gamma}
 \sup_{ x \in X}\sup_{ k
 \le j<i \le 2k}\left |\E \left (\tilde X_i (x)\tilde X_j (x)\right )- \int_X \E\left (\tilde X_i (y)\tilde X_j(y) \right )\nu(dy) \right | <\infty\, .
 \end{equation}
\end{Proposition}

The proof is similar to that of Proposition \ref{prop-aperiodic}, hence is omitted.

\medskip

Finally, we state without proof the following version of Lemma 13 of 
\cite{CDJ}.

\begin{Lemma}\label{startingpoint}
Let $p\ge 2$. Assume that the action is weakly $(\mu,p)$-contracting, that 
 $\mu$ is $F$-adapted  and that 
 $\sigma$ admits a moment of order $p$.
Then
\begin{equation*}
\sup_{ f(x)=f(y)}\left \|\sigma(A_n, x)-\sigma (A_n, y) 
   \right \|_1 = O(1)
   \end{equation*}
for $r \in (1,2]$, 
$$
\sup_{ f(x)=f(y)}\left \|\sigma(A_n, x)-\sigma (A_n, y) 
   \right \|_r  
   =
   \begin{cases}
   O(1)  \quad  \quad \quad \quad \quad \, \text{if $r \leq p-1$}\\
%   O\left( \left(\log n\right )^{1/r} \right)  
% \quad \text{if $ r=p-1 $}\\
   O \left( n^{(r+1-p)/r} \right)
   \quad  \text{if $ r > p-1 $},
  \end{cases}
$$
and for $p \in [2,3]$, 
$$
\sup_{ f(x)=f(y)}\left \|\sigma(A_n, x)-\sigma (A_n, y) 
   \right \|_p  = 
   O\left( n^{1/p} \right)  \, .
$$
\end{Lemma}
\begin{Remark}
If $F$ is trivial, in particular (see Remark \ref{connected}) if $G$ is connected, then the proposition holds with $\sup_{x,y\in X}$ rather than 
$\sup_{f(x)= f(y)}$.
\end{Remark}

%\begin{Definition}
%We shall say that $\sigma$ \emph{does not see $F$} if there exists $(g_1,\ldots , g_{|F|})$  such that $\{s(g_1),\ldots , s(g_{|F|})\}=F$ and 
%for every $g\in G$, every $x\in X$ and every $i\in \{1, \ldots ,|F|\}$ , $\sigma(g,x)=\sigma(g, g_i\cdot x)$.
%\end{Definition}

\subsection{Rates of convergence in the CLT and the ASIP}

Thanks to the results of the previous section we can obtain many probabilistic results, as in \cite{CDJ} or \cite{CDM}. 

Let $p\ge 1$. With the notations and definitions of the previous section, assume that the action of $G$ on $X$ is weakly 
$(\mu,p)$-contracting, that $\mu$ is $F$-strictly aperiodic and that $\sigma$ is a $\R$-valued cocycle admitting a moment of order $p$. Then, the conclusion of Theorems 1 and 2 of \cite{CDJ} hold provided that $S_{n,\overline x}$ is replaced with $\sigma(A_n, x)$.

\medskip

If we consider an $\R^d$-valued cocycle, $d\ge 2$, only the items $(i)$ and $(ii)$ of Theorem 1 of \cite{CDJ} holds without any change. 
 To extend their item (iii) to an $\R^d$-valued cocycle, an application of our Theorem \ref{Thmart} gives the following. 

\begin{Theorem} \label{ThASIPcocycles} Let $p \in ]2,3]$. Assume that the action is weakly $(\mu,p)$-contracting and  that 
$\mu$ is $F$-strictly aperiodic. Let $W_0$ be a random variable, with law $\nu$, independent from $(Y_n)_{n\ge 1}$. Let  
$\sigma$ be an $\R^d$-valued cocycle admitting  a moment of order $p$.  Then the series of matrices
\begin{equation}\label{sigma}
\Sigma = {\rm Var} (\sigma (Y_1, W_0)) + 2 \sum_{k=2}^\infty {\rm Cov}( \sigma (Y_1, W_0), \sigma(Y_k, A_{k-1} W_0))
\end{equation}
converges absolutely. Moreover, 
\begin{enumerate}
\item For   any $\varepsilon >0$, enlarging $\Omega$ if
necessary, there exists a sequence $(N_i)_{i
\geq 1}$ of iid ${\mathbb R}^d$-valued centered gaussian random variables with variance $\Sigma$ such that 
\[
 \sigma (A_n, W_0) - n \lambda_{\mu}- \sum_{i=1}^n N_i  =  \left\{
  \begin{aligned}
  o ( n^{1/p} (\log n)^{ \frac{p+1}{2p} + \varepsilon})& \text{ if $p\in ]2,3[ $}  \\
 o ( n^{1/3} (\log n)^{1 + \varepsilon}) & \text{ if $p=3 $} \\
  \end{aligned}
\right. \text{ almost surely} \, .
\]
\item Enlarging $\Omega$ if
necessary, there exists a sequence $(N _i)_{i\geq 1}$ of iid ${\mathbb R}^d$-valued centered gaussian random variables with variance $\Sigma$ such that 
\[
 \Big \Vert  \sup_{1 \leq k \leq n} \Big |  \sigma (A_k, W_0)  -  k \lambda_{\mu} - \sum_{i=1}^k N_i  \Big |_d   \Big \Vert_1=  \left\{
  \begin{aligned}
  O ( n^{1/p} (\log n)^{ \frac{p-1}{2p} })& \text{ if $p\in ]2,3[ $}  \\
 O ( n^{1/3} (\log n)^{2/3}) & \text{ if $p=3 $} \\
  \end{aligned}
\right. \, .
\]
\end{enumerate}
\end{Theorem}

Moreover an application of our Theorem \ref{ThmartW1cond} (see Remark \ref{remark-ThmartW1cond}) gives the following result concerning rates in the CLT in terms of Wasserstein distance of order $1$. 
\begin{Theorem} \label{ThW1condcocycles} Let $p \in ]2,3]$. Assume that the action is weakly $(\mu,p)$-contracting and that 
$\mu$ is $F$-strictly aperiodic.  Let $W_0$ be a random variable, with law $\nu$, independent from $(Y_n)_{n\ge 1}$. Let $\sigma$ be an $\R^d$-valued cocyle admitting  a moment of order $p$.  Let $\nu_{n}$ be the distribution of $n^{-1/2}  ( \sigma (A_n, W_0)  -  n\lambda_{\mu})$, and let $\Sigma$  be defined by \eqref{sigma}. Then
\[
  W_{1} ( \nu_{n} ,P_{G_{\Sigma}}  \big )  \leq  \left\{
  \begin{aligned}
  &  C  n^{ (2-p) /2} & \text{ if $p \in ]2,3[$} \\
 &  C n^{-1/2}  \log n  & \text{ if $p=3 $} \\
  \end{aligned}
\right.\, , 
\]
where $C$ is a positive constant depending on $(p,d)$ but not on $n$.
\end{Theorem}

\medskip

It is unclear, and probably not true in general, whether the above results hold for $(\sigma(A_n, x))_{n\ge 1}$ (for every $x\in X$) rather than  
 $(\sigma(A_n, W_0))_{n\ge 1}$. However, when  $F$ is the trivial group, using Lemmas \ref{startingpoint-0} or  \ref{startingpoint}, one sees that the above results hold true for  $(\sigma(A_n, x))_{n\ge 1}$, for every $x\in X$ (the same normal variables being used for all $x$). 

%\medskip

%Another situation where the results hold for $(\sigma(A_n, x))_{n\ge 1}$, for every $x\in X$ is under the following assumption: for every $\f_1,\f_2\in F$ 
%and every $x\in X_{\f_1}$, there exists $y\in X_{\f_2}$ such that $\sigma(g,x)=\sigma(g,y)$ for every $g\in G$. This hypothesis will be verified in our applications, namely for the Iwasawa cocycle. 

\medskip

\subsection{Applications to the Iwasawa cocycle for reductive Lie groups} \label{Iwasawa-section}

In this section, we shall give a general situation to which the previous 
sections apply. 

\medskip

Our presentation will borrow notations as well as results from the monograph by 
Benoist and Quint \cite{BQ-book}. We refer  to \cite{BQ-book} for 
any complements on the topic discussed in that section. However, for the sake of clarity, we shall not cover the generality treated in 
\cite{BQ-book}. In particular, we shall only condider algebraic groups over the field of real numbers, while the results presented here 
extend to local fields modulo several technical issues. 

\medskip

We give the definition of the Iwasawa cocycle in a general setting and we will make explicit all the objects in the particular case of the $d$-dimensional linear group.

\medskip

Let $G$ be a reductive algebraic real Lie group with lie algebra $\g$. 
Recall that $G$ is said to be reductive if its unipotent radical, 
that is the greatest connected normal subgroup of $G$ whose elements are 
unipotent, is reduced to $\{e\}$, where $e$ stands for the neutral element of $G$. Recall also that $G$ is algebraic if it is the set of solutions of a (finite) system of polynomial equations (over $\R$).

Define the killing form on $\g\times \g$ by 
$$
 {\rm killing}(x,y):={\rm tr} ({\rm ad}\, x\, {\rm ad}\, y)\qquad \forall x,y\in \g\, .
$$

Let 
$K$ be a maximal compact subgroup of $G$ with Lie algebra $\k$. Let $\s$ be the orthogonal space of $\k$ for the killing form and let 
$\mfa$ be a Cartan subspace of $\s$, that is $\mfa$ is a commutative 
 subalgebra 
of $\s$ whose elements are diagonalizable over $\R$ and which is maximal with those properties. Let $A:={\rm exp}\, \mfa$ be the corresponding connected lie subgroup of $G$.

Let $U$ be a maximal unipotent subgroup of $G$ that is normalized by $A$. 
Let $P:=N_G(U)$ be the normalizer of $U$ in $G$.
% and let $Z:=Z_G(A)$ be the centralizer of $A$ in $G$.

\medskip

Denote by $G_c$, $K_c$ and $P_c$  the respective  connected 
components of  $G$, $K$ and  $P$.

\smallskip

We have the Iwasawa decomposition (the proof is sketched page 130 of 
\cite{BQ-book})

\begin{equation}\label{iwasawa-dec}
G=KP_c=KA U\, .
\end{equation}

The compact space $\PP:= G/P_c$ is known as the flag variety of $G$. 
%The decomposition \eqref{iwasawa-dec} is not unique since we may have 
%$K\cap Z_c\neq\{e\}$. However, g
Given $\eta=kP_c$ with $k\in K$ and $g\in G$,  there exists a unique element $z$ of $\mfa$ such that $gk\in K{\rm exp}\, z\, U$.  We denote by $\sigma(g,\eta)$  that unique element. 

Then, $\sigma$ is a continuous cocycle known as the Iwasawa cocycle. It is related to the Cartan projection that we shall introduce now. 

\medskip

Let $\Sigma$ be the set of roots associated with $\mfa$, that is, denoting by $\mfa^*$ the dual vector space of $\mfa$,

$$\Sigma:=\{\alpha \in \mfa^*-\{0\}\,:\, \g^\alpha\neq \{0\}\}\, ,
$$
where 
$$
\g^\alpha :=\{y\in \g \,:\, \forall \, x \in \mfa,\, {\rm ad}\, x (y)\, 
=\alpha(x)y 
\}\, .
$$
Let $\Sigma^+$ be the set of positive roots associated with $U$, that is, 
if $\u$ stands for the lie algebra of $U$, $\Sigma^+$ is characterized by
$$
\u=\oplus_{\alpha \in \Sigma^+}\g^\alpha\, .
$$

Define then $\mfa^+:=\{x\in \mfa\, :\,\forall \alpha\in \Sigma^+,\, \alpha(x)\ge 0\}$. We have the Cartan decomposition (see also page 130 of 
\cite{BQ-book})

$$
G=K {\rm exp}\, \mfa^+\, K_c\, .
$$

Finally, define the Cartan projection as follows: for every $g\in G$, 
let $\kappa(g)$ be the unique $z\in \mfa^+$ such that $g\in K \,
{\rm exp}\, z\, K_c$.

\medskip

For people not familiar with Lie groups we consider now the case where $G={\rm Gl}_d(\R)$, $d\ge 2$, and describe the different spaces introduced above. 

\medskip

 In that case, one may take $K={\rm O}_d(\R)$ the orthogonal group, $A$ the group of invertible diagonal matrices of size $d$ with positive entries, $U$ the group of  upper triangular matrices whose diagonal terms are all 1 and $P$ the subgroup of invertible upper triangular matrices.

Let us describe the corresponding Lie algebras are. We have $\g =M_d(\R)$, the set of  matrices of size $d$.  Then,  the Killing form is given by ${\rm killing}(M,N) =2d \, {\rm tr}(MN)-2{\rm tr} M{\rm tr} N$ for every 
$M,N\in \g $.   Then,  $\k= \{M\in \g\, :\, M=-M^t\}$  the set of antisymmetric matrices; 
$\s$ is the set of symmetric matrices of size $d$  and 
$\mfa $ the set of diagonal real matrices of size $d$.   Then, 
$G_c=\{M\in G\, :\, {\rm det}\, M>0\}$, $K_c={\rm SO}_d(\R)$ and $P_c$ is the group of upper triangular matrices with non-negative entries on 
the diagonal. 

\medskip

Finally, the set of positive roots is given by $\Sigma^+=\{\alpha_{ij}\, :\, 1\le i<j\le d\}$ where for $M={\rm diag}(a_m)_{1\le m\le d}$, $\alpha_{ij}(M)
=a_i-a_j$ and $\u$ is the set of upper triangular matrices with vanishing diagonal. In particular,   $\mfa^+$ consists of the set of diagonal matrices with entries in non-increasing order.

\medskip

In this case, the Iwasawa decomposition may be seen as a by-product the Gram-Schmidt orthonormalization algorithm. Indeed, if $M\in G$ then 
$M$ sends the canonical basis to a non necessarily orthonormal basis. Orthonormalizing the latter thanks to the Gram-Schmidt algorithm one sees that there exists an upper triangular matrix $T$ with positive diagonal coefficients such that $MT$ sends the canonical basis to an orthonormal basis, i.e. 
$MT\in O_d(\R)$.

\medskip

The Cartan projection also has a nice interpretation in this case. It may be seen as a by-product of the polar decomposition. In particular, for every $M\in G$ there exist unique matrices $K,K'\in O_d(\R)$ and $A\in G$ such that ${\rm det}\, K'=1$, $A$ is diagonal with positive diagonal coefficients in non 
decreasing order. The coefficients $(\lambda_i)_{1\le i\le d}$ of $A$ are the square roots of the eigenvalues of $M^tM$. Hence, 
$\kappa(M)={\rm diag }\, (\log \lambda_i)_{1\le i\le d}$.

\medskip

\medskip

Let us come back to the general case. It follows from section 13.1 of \cite{BQ-book} that $\PP$ may be endowed 
with a metric $d$  compatible with the quotient topology on $\PP$. The 
metric is defined by (13.3) (notice that in our situation, i.e. 
$G$ is a real algebraic Lie group, $\Theta=\Pi$, where 
$\Pi$ is the set of simple roots, that is the roots of $\Sigma^+$ that 
are not the sum of two roots of $\Sigma^+$). Moreover, 
$\PP$ is fibered over $F$ where $F:=G/G_c$ (see page 142 of 
\cite{BQ-book}. The group $F$ is finite by Lemma 6.21 of \cite{BQ-book}.  

\medskip

It follows from (13.4) of \cite{BQ-book} that there exists 
$C,D>0$ such that 

$$
d(g\cdot x,g\cdot y)\le C{\rm e}^{D|\kappa(g)|_\mfa}d(x,y) \qquad \forall g\in G,\, \forall x,y\in \PP\,,
$$
where $|\cdot |_\mfa$ is an euclidean norm on $\mfa$. In particular, for every $g\in G$, 
$$
{\rm Lip}(g)\le C{\rm e}^{D|\kappa(g)|_\mfa}\, .
$$

\begin{Definition}
Following \cite[Section 9]{BQ-book}, we say that a Borel probability measure on $G$ is Zariski dense if 
the subsemigroup  $\Gamma_\mu$ spanned by the support of $\mu$ is 
$G$.
\end{Definition}
If $\mu$ is a  Zariski dense Borel probability measure on $G$, 
then  (13.10) of \cite{BQ-book} holds, which means exactly that 
\eqref{contracting} is satisfied. In particular, for a Zariski dense 
Borel probability measure $\mu$ on $G$, the actions of $G$ on $\PP$
 is weakly $(\mu,p)$-contracting for some $p\ge 1$, as soon as 
 \begin{equation}\label{pth-moment}
 \int_G|\kappa(g)|^p_\mfa \, \mu(dg)<\infty\, .
 \end{equation}

In order to apply our results, it remains to control the Iwasawa 
cocycle.  It follows from (8.16) and (13.5) of 
\cite{BQ-book}, that there exist $C,D>0$ such that 
\begin{gather*}
|\sigma(g,\eta)|_\mfa \le C|\kappa(g)|_\mfa \qquad \forall g\in G,\, 
\eta \in \PP\\
\sigma_{\rm Lip}(g)\le C{\rm e}^{D|\kappa(g)|_\mfa}\qquad \forall 
g\in G\, .
\end{gather*}

In particular, $\sigma$ admits a moment of order $p$, for some $p\ge 1$ as soon as \eqref{pth-moment} holds and our results apply.

\begin{Theorem}\label{ThASIP-Iwasawa}
Let $G$ be a reductive algebraic Lie group. Let $\mu$ be a Zariski dense Borel probability measure on $G$ that is $F$-strictly aperiodic. Let $\sigma$ be the Iwasawa cocycle and 
$\kappa$ be the Cartan projection as above. Assume that \eqref{pth-moment} holds for some $p\in (2,3]$, and let $\Sigma$ be defined by  \eqref{sigma}. Then, 
\begin{enumerate}
\item For any $\varepsilon >0$, enlarging $\Omega$ if
necessary, there exists a sequence $(N_i)_{i
\geq 1}$ of iid $\mfa$-valued centered gaussian random variables with variance $\Sigma$  such that 
\[
 \big|V_n - \sum_{i=1}^n N_i  \big|_\mfa =  \left\{
  \begin{aligned}
  o ( n^{1/p} (\log n)^{ \frac{p+1}{2p} + \varepsilon})& \text{ if $p\in ]2,3[ $}  \\
 o ( n^{1/3} (\log n)^{1 + \varepsilon}) & \text{ if $p=3 $} \\
  \end{aligned}
\right. \text{ almost surely} \, ,
\]
where $(V_n)_{n\ge 1}$ is either of the processes $(\sigma(A_n, W_0)-n\lambda_\mu)_{n\ge 1})$, $(\kappa(A_n)-n\lambda_\mu)_{n\ge 1}$ or $(\sigma(A_n, x)-n\lambda_\mu)_{n\ge 1}$ for a given 
$x\in {\mathcal P} $.
\item Enlarging $\Omega$ if
necessary, there exists a sequence $(N _i)_{i
\geq 1}$ of iid ${\mathbb R}^d$-valued centered gaussian random variables with variance $\Sigma$ such that 
\[
 \Big \Vert  \sup_{1 \leq k \leq n} \big |   V_k - \sum_{i=1}^k N_i  \big |_\mfa   \Big \Vert_1=  \left\{
  \begin{aligned}
  O ( n^{1/p} (\log n)^{ \frac{p-1}{2p} })& \text{ if $p\in ]2,3[ $}  \\
 O ( n^{1/3} (\log n)^{2/3}) & \text{ if $p=3 $} \\
  \end{aligned}
\right. \, ,
\]
where $(V_n)_{n\ge 1}$ is either of the processes $(\sigma(A_n, W_0)-n\lambda_\mu)_{n\ge 1})$ or $(\kappa(A_n)-n\lambda_\mu)_{n\ge 1}$. 
If moreover $G$ is connected or $G=GL_d(\R)$, then we also have 
\[
\sup_{x\in X} \Big \Vert  \sup_{1 \leq k \leq n} \big |  \sigma (A_k, x)  -  k \lambda_{\mu} - \sum_{i=1}^k N_i  \big |_\mfa   \Big \Vert_1=  \left\{
  \begin{aligned}
  O ( n^{1/p} (\log n)^{ \frac{p-1}{2p} })& \text{ if $p\in ]2,3[ $}  \\
 O ( n^{1/3} (\log n)^{2/3}) & \text{ if $p=3 $} \\
  \end{aligned}
\right. \, .
\]
\end{enumerate}
\end{Theorem}
\begin{Remark}
The almost sure set in Item $1.$ above depends on $x$ when $V_n=\sigma(A_n, x)$.
\end{Remark} 

\begin{Theorem} \label{ThW1-Iwasawa} Let $G$ be a reductive algebraic Lie group. Let $\mu$ be a Zariski dense Borel probability measure on $G$ that is 
$F$-strictly aperiodic. Let $\sigma$ be the Iwasawa cocycle and 
$\kappa$ be the Cartan projection as above. Assume that \eqref{pth-moment} holds for some $p\in (2,3]$, and let $\Sigma$ be defined by  \eqref{sigma}. Let $G_{\Sigma}$ be an $\mfa$-valued centered gaussian random variables with variance $\Sigma$. Then 
\[
 W_{1} (P_{n^{-1/2} V_n} , P_{G_{\Sigma}}  \big )  \leq  \left\{
  \begin{aligned}
  &  C  n^{ (2-p) /2} & \text{ if $p \in ]2,3[$} \\
 &  C n^{-1/2}  \log n  & \text{ if $p=3 $} \\
  \end{aligned}
\right.\, , 
\]
where $C$ is a positive constant depending on $(p,d)$ but not on $n$ and $(V_n)_{n\ge 1}$ is either of the processes $(\sigma(A_n, W_0)-n\lambda_\mu)_{n\ge 1})$, $(\kappa(A_n)-n\lambda_\mu)_{n\ge 1}$. 
If moreover $G$ is connected or $G=GL_d(\R)$, then we also have 
\[
\sup_{x\in X}   W_{1} (P_{n^{-1/2} (\sigma(A_n, x)-n\lambda_\mu)} , P_{G_{\Sigma}}  \big )  \leq  \left\{
  \begin{aligned}
  &  C  n^{ (2-p) /2} & \text{ if $p \in ]2,3[$} \\
 &  C n^{-1/2}  \log n  & \text{ if $p=3 $} \\
  \end{aligned}
\right.\, , 
\]
\end{Theorem}

\section{Proofs}

\setcounter{equation}{0}

All along the proofs we denote by $C$ a numerical constant  which may vary
from line to line and which may depend on $d$ and $p$ but not on $n$.  We shall also denote sometimes by $\E_i$ the conditional expectation with respect to $\F_i$. 

\subsection{Proofs of the results of Section \ref{Sec:mart}}

\subsubsection{Preliminaries} \label{sectpreli}

Let $\Sigma= \E (d_0 d_0^t)$. Suppose that $\Sigma$ is nonnull (otherwise there is nothing to prove since the $d_k$'s are all almost surely equal to zero). Since 
$\Sigma$ is symmetric and positive-semidefinite, it follows that there exists  a $d$-dimensional  orthogonal matrix $P$ such that 
\[
\Sigma = P D P^t \, , 
\]
with $D = {\rm Diag} ( \lambda_1, \ldots, \lambda_d)$ where the $\lambda_i$'s are the eigenvalues of $\Sigma$ ranking in the non-increasing order. All these eigenvalues are reals and non-negative. Let $m  \in \{1, \ldots, d \}$  be the number of eigenvalues that are positive and $\Delta = {\rm Diag} ( \lambda_1, \ldots, \lambda_m, \lambda_m, \ldots, \lambda_m)$ (i.e. the $d$-dimensional  diagonal matrix such that the first $m$ diagonal elements are equal to the first $m$ diagonal elements of $D$ and the others to $\lambda_m$). Denote by ${\rm J}_m$  the $d$-dimensional
diagonal matrix such that the first $m$ diagonal elements are equal to $1$ and the others to $0$. Since $D= \Delta^{1/2} {\rm J}_m \Delta^{1/2}$, 
\[
\Sigma = P  \Delta^{1/2} {\rm J}_m \Delta^{1/2} P^t \, . 
\]
Setting $\Gamma =  {\rm Diag} ( \lambda^{-1/2}_1, \ldots,  \lambda^{-1/2}_m,  \lambda^{-1/2}_m, \ldots,  \lambda^{-1/2}_m)$ (i.e. the inverse of $\Delta^{1/2}$) and $A =\Gamma P^t $, it follows that \[
A \Sigma A^t = {\rm J}_m \, , 
\]
since $PP^t = P^t P =  {\rm I}_d$ where as usual $ {\rm I}_d$ denotes the identity matrix on ${\mathbb R}^d$. Note that $A$ is invertible and $A^{-1} = P  \Delta^{1/2}$. 
For any integer $k$, let now
\[
m_k = A d_k \, .
\]
Note that $(m_n)_{n \in {\mathbb N}}$ is a ${\mathbb R}^d$-valued stationary sequence of  martingale differences   with respect to $(\F_n)_{n \in {\mathbb N}}$ such that $\bkE |m_0 |_d^p < \infty$ and satisfying 
$ \E (m_0 m_0^t)= {\rm J}_m$. Hence, clearly, $(m_0)_i =0$ a.s. for any $i=m+1, \ldots, d$. For any integer $k$, let now $d'_k$ be the ${\mathbb R}^m$-valued random vector whose components are equal to the $m$-first components of $m_k$, that is
\[
d_k' = \big ( (m_k)_1, \ldots,  (m_k)_m\big )^t \, .
\]
Clearly, $(d'_n)_{n \in {\mathbb N}}$ is a ${\mathbb R}^m$-valued stationary sequence of  martingale differences   with respect to $(\F_n)_{n \in {\mathbb N}}$ such that $\bkE |d'_0 |_m^p < \infty$ and satisfying 
$ \E (d'_0 (d'_0)^t)= {\rm I}_m$. Let $M'_n = \sum_{i=1}^n d'_i$. 

A common key result for the proofs of Theorems \ref{Thmart} and \ref{ThmartW1cond} is the following lemma:  Let  ${\rm Lip}(|\cdot |_m, {\mathcal F}_{0})$ be the set of measurable functions $g:{\mathbb R}^{m}
\times  \Omega \rightarrow {\mathbb R}$ with respect to the $\sigma$-fields 
${\mathcal B} ( {\mathbb R}^{m}) \otimes {\mathcal F}_{0} $ and ${\mathcal B} ({\mathbb R})$, 
such that $g( \cdot, \omega) \in {\rm Lip}( |\cdot |_{m})$ and $g(0,\omega)=0$  for any 
$\omega \in \Omega$. 
For the sake of brevity, we shall write $g(x)$ in place of $g(x, \omega)$. 
\begin{Lemma} \label{keylemma} Under the assumptions of Theorem \ref{Thmart}, setting $T_n = \sum_{i=1}^n N_i$ where $(N_i)_{i
\geq 1}$ is a sequence of iid ${\mathbb R}^m$-valued centered gaussian random variables with ${\rm Var} (N_1) = {\rm I}_m$, we have
\[
\sup_{g \in {\rm Lip}( |\cdot |_{m}, {\mathcal F}_{0})} {\mathbb E}  (  g(M'_n)  ) 
- {\mathbb E} (g( T_n))  \leq   \left\{
  \begin{aligned}
  C  n^{ (3-p) /2} & \text{ if $p \in ]2,3[$} \\
  C \log n  & \text{ if $p=3 $} \\
  \end{aligned}
\right.\, , 
\]
where  $C$ is a positive constant depending on $(p,m)$ but not on $n$.
\end{Lemma}

\subsubsection{Proof of Theorem \ref{Thmart}} 

The construction of the approximating sequence of Gaussian random variables uses the ideas developed in the proof of Theorem 2.1 in Merlev\`ede and Rio \cite{MR}.

\medskip

Let $(m(L))_{L \in {\mathbb N}} $ be a  sequence of non-negative integers that will be specified later but such that $m(L)\leq L$. Let 
\begin{equation*} \label{defUkL}
I_{k,L} = ]2^L + (k-1)2^{m(L)} ,  2^L + k 2^{m(L)}] \cap {\mathbb N} \ \text{and}\ U_{k,L} = \sum_{i\in I_{k,L}} d'_i \, ,
\, k \in \{1, \cdots, 2^{L-m(L)} \} \, .
\end{equation*}
Let $P_{U_{k,L} | {\mathcal F}_{2^L + (k-1)2^{m(L)}
}}$ be the conditional law of $U_{k,L}$ given $\mathcal{F}_{2^L +
(k-1)2^{m(L)} }$ and ${\mathcal N}_{2^{m(L)}}$ denote the ${\mathcal N} ( 0, 2^{m(L)} {\rm I}_m)$-law. The probability space is assumed to be large enough to contain a sequence $(\delta_i)_{i \in {\mathbb Z}}$ of iid random variables uniformly distributed on $[0,1]$,  independent of the sequence $(d_i)_{ i \in {\mathbb Z}}$ (otherwise we enlarge it). According to R\"uschendorf \cite{Rus} (see also Theorem 2 in \cite{DPR}),  
there exists a ${\mathbb R}^m$-valued random variable $V_{k,L}$ with law ${\mathcal N}_{2^{m(L)}}$, measurable with respect to $\sigma(\delta_{2^L +k 2^{m(L)}}) \vee \sigma(U_{k ,L}) \vee {\mathcal F}_{2^{L}+ (k-1)2^{m(L)}} $, independent of ${\mathcal F}_{2^L + (k-1)2^{m(L)}}$ and such that
\begin{eqnarray} \label{coupling1}
{\mathbb E} \big (  \big | U_{k,L} - V_{k,L} \big |_m \big ) & = & {\mathbb E} \big ( W_{1} (P_{U_{k,L} | {\mathcal F}_{2^L + (k-1)2^{m(L)}
}} ,  {\mathcal N}_{2^{m(L)}} ) \big ) \\
& = &  \bkE  \sup_{f \in {\rm Lip} (|\cdot |_{m} ) } \Big (  {\mathbb E} \big (  f(U_{k,L}) | {{\mathcal F}_{2^L + (k-1)2^{m(L)}}} \big ) - {\mathbb E} (f(V_{ k,L})) \Big ) \nonumber \, ,
\end{eqnarray}
where we recall that $ {\rm Lip} (| \cdot |_m)$ is the set of functions from ${\mathbb R}^m$ into ${\mathbb R}$ that are $1$-Lipschitz with respect to the euclidian norm $|\cdot|_{m} $ on ${\mathbb R}^m$.  From Point 2 of Theorem 1 in \cite{DPR}, the following inequality  holds: 
\[
 {\mathbb E} \big ( W_{1} (P_{U_{k,L} | {\mathcal F}_{2^L + (k-1)2^{m(L)}
}} ,  {\mathcal N}_{2^{m(L)}} ) \big )  
= \sup_{g \in {\rm Lip}(|\cdot|_{m}, {\mathcal F}_{2^L})} {\mathbb E}  (  g(U_{1,L} )  ) 
- {\mathbb E} (g( V_{1,L}  )) \, .
\]
Hence, using Lemma \ref{keylemma},  we have that
for any $L \in {\mathbb N}$ and any $ k \in \{1, \cdots, 2^{L-m(L)} \}$,
\beq \label{cons1keylemma}
{\mathbb E} \big (  \big | U_{k,L} - V_{k,L} \big |_m \big )  \leq   \left\{
  \begin{array}{rcr}
   C  2^{  (3-p) m(L) /2} & \text{ if $p \in ]2,3[$} \\
  C  m(L) & \text{ if $p=3 $} \\
  \end{array}
\right.\, , 
\eeq
where $C$ is a positive constant depending on $(p,m)$ but not on $(k,L)$.

By induction on $k$, the random variables $(V_{k,L})_{k= 1, \ldots, 2^{L-m(L)}}$
are mutually independent, independent of ${\mathcal F}_{2^L}$ and with law ${\mathcal N}_{2^{m(L)}}$ . Hence we have constructed Gaussian random variables $(V_{k,L})_{L \in {\mathbb N}, k=1, \ldots, 2^{L-m(L)}}$ that are mutually independent.
Now we construct a
sequence $(Z_i)_{i \geq 1}$ of iid standard Gaussian random vectors in ${\mathbb R}^m$. For any
$L \in {\mathbb N}$ and any $k \in \{1, \cdots, 2^{L -m(L)} \}$ 
the random variables \[(Z_{2^L + (k-1)2^{m(L)} +1}, \ldots , Z_{2^L + k 2^{m(L)}}) \]
are defined in the following way. If $m(L)=0$, then 
$Z_{2^L + k 2^{m(L)}} = V_{k,L}$.
If $m(L)>0$, then by the Skorohod lemma \cite{Sko}, there exists a measurable function $g$ from ${\mathbb R}^m \times [0,1]$ 
in $({\mathbb R}^{m})^{\otimes 2^{m(L)}}$ such that, for any pair 
$(V,\delta)$ of independent random variables with respective laws  ${\mathcal N}_{2^{m(L)}}$ and 
the uniform distribution over $[0,1]$, $g(V,\delta) = (N^t_1, \ldots , N^t_{2^{m(L)}})$ is a Gaussian random vector with iid components such that $V = N_1 + \cdots + N_{2^{m(L)}}$ a.s. 
Next we set
$$
(Z^t_{2^L + (k-1)2^{m(L)} +1}, \ldots , Z^t_{2^L + k 2^{m(L)}})  = g (V_{k,L} , \delta_{2^L + (k-1)2^{m(L)} +1} ) \, . 
$$
We have then constructed a sequence $(Z_i)_{i \geq 2}$ of iid standard Gaussian random vectors in ${\mathbb R}^m$ such that, for any $L \in {\mathbb N}$ and any $ k \in \{1, \cdots, 2^{L-m(L)} \}$,
\[
 V_{k,L} =  \sum_{i\in I_{k,L}}  Z_i  \  \ a.s. 
\]
To complete the construction of the sequence  $(Z_i)_{i \geq 1}$, it suffices to consider a ${\mathbb R}^m$-valued standard Gaussian random vector $Z_1$ independent of $(d_i, \delta_i)_{i \in {\mathbb Z}}$ which is always possible by enlarging enough the underlying probability space. 

Let us now complete the proof of Theorem \ref{Thmart}.  First, for any $k \geq 1$, we set 
\[
Y_k = A^{-1} Z_k' \mbox{ where } Z_k' = ( Z_k^t, {\bf 0}_{d-m} )^t \, , 
\]
${\bf 0}_{d-m}$ denoting the row vector of dimension $d-m$ whose all components are equal to $0$. Since ${\rm Var} (Z_k') = {\rm J}_m$, we get that 
\[
{\rm Var} (Y_k) = A^{-1} {\rm J}_m(A^{-1})^t = \Sigma \, .
\]
So,  $(Y_i)_{i \geq 1}$ is a sequence of iid centered Gaussian random vectors in ${\mathbb R}^d$ with covariance matrix $\Sigma$. Note that 
\begin{multline*}
\sup_{k \leq n} \big |  M_k - \sum_{i=1}^k Y_i \big |_d = \sup_{k \leq n} \big |  A^{-1}  \big ( \sum_{i=1}^k m_i - \sum_{i=1}^k Z'_i  \big ) \big |_d  \\
 \leq  \Vert A^{-1}  \Vert_2  \sup_{k \leq n} \big |  \sum_{i=1}^k m_i - \sum_{i=1}^k Z'_i  \big |_d  =   \sqrt{\lambda_1}
\sup_{k \leq n} \big |  \sum_{i=1}^k d'_i - \sum_{i=1}^k Z_i  \big |_m  \, , 
\end{multline*}
where for the last inequality we have used the fact that 
\[
\Vert A^{-1}  \Vert_2 = \sqrt{\rho ((A^{-1})^t A^{-1}) } = \sqrt{ \rho (   \Delta^{1/2} P^t P  \Delta^{1/2}) }  =  \sqrt{ \rho (   \Delta) } = \sqrt{\lambda_1} \, .
\]
Above and in the rest of the paper, for any $B \in {\mathcal M}_n ( {\mathbb C})$, the notation $\rho ( B)$ means the spectral radius of $B$.   To prove the theorem, it suffices then to show that 
\beq \label{but1pthm}
\sup_{k \leq n} \big |  \sum_{i=1}^k d'_i - \sum_{i=1}^k Z_i  \big |_m =  o ( n^{1/p} (\log n)^{a}), \text{ almost surely} \, ,
\eeq
for a suitable $a$. 
With this aim, set 
$S_j = \sum_{i=1}^j d'_i$ and $T_j = \sum_{i=1}^j Z_i$ and  let 
\begin{eqnarray} \label{0dec} 
D_L := \sup_{ \ell \leq 2^{L}} | \sum_{i = 2^L +1}^{2^L + \ell} (d'_i -Z_i)|_m \, .
\end{eqnarray}
Let $N\in {\BBN}^*$ and  $k  \in ]1, 2^{N+1}]$.
We first notice that $D_L \geq | (S_{2^{L + 1} } -  T_{2^{L + 1} }) - ( S_{2^L } - T_{2^L } )|_m$,
so that, if $K$ is the  integer such that $2^K < k \leq  2^{K+1}$, $|S_k - T_k|_m \leq |d'_1 - Z_1|_m+ D_0 + D_1 + \cdots + D_K$. 
Consequently since $K \leq N$,
\begin{eqnarray} \label{1dec} \sup_{1 \leq k \leq 2^{N+1} }|S_k -T_k|_m  \leq 
|d'_1 - Z_1|_m + D_0 + D_1 + \cdots + D_N \, . 
\end{eqnarray}
In addition,  the following decomposition is valid: 
\beq \label{decsup} 
D_L  \leq D_{L,1} + D_{L,2} \, ,
 \eeq 
where 
\[
D_{L,1}:= \sup_{k \leq 2^{L-m(L)} } \Big|
\sum_{ \ell =1}^k (U_{\ell,L} - V_{\ell,L}) \Big|_m
\ \text{and}\ 
D_{L,2}:= \sup_{ k \leq 2^{L-m(L)}} \sup_{ \ell \in I_{k,L} } 
\Big| \sum_{i = \inf  I_{k,L}  }^{\ell} (d'_i -Z_i) \Big|_m \, .
\]
\noindent{\bf End of the proof of Item 1.} 
Let $\varepsilon >0$ and 
\[
\text{$a_p=\frac{1}{2} +\frac{1}{2p} + \varepsilon$  if $p \in ]2, 3[$  and $a_3= 1 + \varepsilon$.} 
\] 
From \eqref{0dec} and \eqref{1dec}, it follows  that the proof of Item 1 will be complete if we can show that, for any $L \in {\mathbb N}$,
\beq \label{but2pthm}D_{L,1}=  O (2^{L/p} L^{a_p} ) \ \, \text{  and } \ \, D_{L,2}=  O (2^{L/p} L^{a_p} ) \ \, \text{ a.s.}
\eeq
As we shall see below, this will be achieved by selecting the sequence $(m(L))_{L \geq 0}$ as follows in the construction of the iid gaussian vectors $(Y_i)_{i \geq 1}$ as described above: set 
$b_p = \frac{1}{p} $  if $p \in ]2, 3[$  and $b_3= 1 $ and 
\beq \label{defml}
m(L) = \Big [ \frac{2L}{p} + b_p \log_2 L \Big ] \, ,\ \text{so that}\ 
\frac 12 2^{2L/p}L^{b_p}  \leq 2^{m(L)} \leq 2^{2L/p} L^{b_p} \, ,
\eeq 
square brackets designating as usual the integer part and $\log_2 (x)= (\log x)/(\log 2)$.

\smallskip

To prove the first part of \eqref{but2pthm}, note that, by \eqref{cons1keylemma} and the selection of $m(L)$, we get that, for any $c >0$, 
\begin{multline*}
{\mathbb P} \big ( D_{L,1} \geq c 2^{L/p} L^{a_p}  \big )    \leq  \frac{c}{ 2^{L/p} L^{a_p}}  
\sum_{\ell=1}^{2^{L-m(L)}}  {\mathbb E} \big ( \big | U_{\ell,L} - V_{\ell,L} \big|_{m}  \big )  \leq C \frac{2^L}{2^{L/p} L^{a_p}}  2^{ (1-p)m(L)/2}  ( 1 + {\bf 1}_{p=3} m(L) ) \\
 \leq  \frac{C }{L^{a_p +(p-1)b_p/2 }}   ( 1 + {\bf 1}_{p=3} L  ) \, . 
\end{multline*}
Hence, for any $c >0$, 
\[
{\mathbb P} \big ( D_{L,1} \geq c 2^{L/p} L^{a_p}  \big )      \leq  \frac{C }{L^{ 1 + \varepsilon }}  \, , \]
which together with the Borel-Cantelli lemma implies the first part of \eqref{but2pthm}.

We turn now to the proof of the second part of \eqref{but2pthm}. With this aim, we set
\[
x_L = \kappa 2^{L/p} L^{a_p}  \text{ for some $\kappa >0$}\, , 
\]
and first notice that, by stationarity, for any $y >0$,
\beq \label{decdl2} {\mathbb P}
( D_{L,2} \geq 2 y) \leq 2^{L-m(L)} {\mathbb P} \Bigl(  \sup_{
\ell \leq 2^{m(L)}} | S_{\ell} |_m \geq y \Bigr) +  2^{L-m(L)} {\mathbb P} 
\Bigl(  \sup_{\ell \leq 2^{m(L)}} | T_{\ell} |_m \geq y
\Bigr)  \, .\eeq
By L\'evy's inequality (see for instance Proposition 2.3 in \cite{LeTa}),
\begin{equation} {\mathbb P} \Bigl(  \sup_{
 \ell\leq 2^{m(L)}} | T_{\ell} |_m \geq x_L
\Bigr) \leq  \sum_{i=1}^m  {\mathbb P} \Bigl(  \sup_{
 \ell\leq 2^{m(L)}} | (T_{\ell})_i |  \geq m^{-1/2} x_L
\Bigr)  \leq 2 m  \exp  \Bigl( - \frac{x_L^2}{ m 2^{ m(L) +1}
} \Bigr)   \, . \label{LI}
\end{equation}
On the other hand,  using Proposition \ref{inegamax} of the appendix, we get that, for any integer $i \in [1, d ]$,  there exist two positive constants $c_1$ and $c_2$ such that,  for any $x >0$, 
\beq \label{inegamaxmart}
\BBP  \Bigl(  \sup_{ \ell \leq n } | (M_{\ell})_i | \geq x \Bigr)  \leq  
c_1\exp \Bigl( - \frac{x^2}{c_2 n } \Bigr) + c_2 n x^{-p}   \, .
\eeq
Now note that 
\[
 \sup_{
\ell \leq 2^{m(L)}} | S_{\ell} |_m =  \sup_{
\ell \leq 2^{m(L)}} | A M_{\ell} |_d \leq \Vert A \Vert_2  \sup_{
\ell \leq 2^{m(L)}} |  M_{\ell} |_d  \, .
\]
But $\Vert A \Vert_2  =  \sqrt{\rho(P \Gamma^2 P^t) }   = \lambda_m^{-1/2} $. Therefore, applying  inequality \eqref{inegamaxmart}, we get 
\begin{multline} \label{consinegamaxmart} 
{\mathbb P} \Bigl(  \sup_{
\ell \leq 2^{m(L)}} | S_{\ell} |_m \geq x_L \Bigr)  \leq \sum_{i=1}^m  {\mathbb P} \Bigl(  \sup_{
 \ell\leq 2^{m(L)}} | (M_{\ell})_i |  \geq ( m^{-1/2} \lambda_m^{1/2} x_L  \Bigr)  
\Bigr) \\
 \leq   c_1 m \exp \Bigl( - \frac{x_L^2 \lambda_m }{c_2  m   2^{m(L) }} \Bigr) +  c_2 m^{(2+p)/2} \lambda_m^{-p/2} 2^{m(L)} x_L^{-p} \, .
\end{multline}
Starting from \eqref{decdl2} and considering the upper bounds \eqref{LI} and \eqref{consinegamaxmart}, it follows that 
\beq \label{decdl2bis} {\mathbb P}
( D_{L,2} \geq 2 x_L) \leq m ( c_1 + 2 )  2^{L-m(L)}\exp \Bigl( - \frac{  x_L^2 \lambda_m }{   c(m)   2^{m(L) }} \Bigr)  +  c_2 m^{(2+p)/2} \lambda_m^{-p/2} 2^{L} x_L^{-p}  \, ,\eeq
where 
\[
c(m) = m  \max ( 2, c_2 \lambda_m^{-1}  )\, .
\] 
For any choice of $\kappa$, by the selection of $x_L$ and  since $p a_p >1$, it follows that $ \sum_{L \geq 1} 2^{L} x_L^{-p}  < \infty$.  On another hand,
\[
\frac{  x_L^2  }{    2^{m(L) }} \geq \kappa^2 L^{ 2 a_p -b_p}  = \kappa^2 L^{ 1+  2 \varepsilon } \, .
\]
It follows that for any $\kappa >0$,
\[
 \sum_{L \geq 1} 2^{L} \exp \Bigl( - \frac{  x_L^2 \lambda_m }{   c(m)   2^{m(L) }} \Bigr)   < \infty \, .
\]
So, overall, starting from \eqref{decdl2bis} and using the Borel-Cantelli lemma, we can conclude that  the second part of \eqref{but2pthm} holds. This ends the proof of Item 1.

\medskip

\noindent{\bf End of the proof of Item 2.} 
Let $\varepsilon >0$ and 
\[
\text{$a_p=\frac{1}{2} - \frac{1}{2p}$  if $p \in ]2, 3[$  and $a_3= 2/3$.} 
\] 
Starting again from, \eqref{0dec} and \eqref{1dec}, it follows  that the proof of Item 2 will be complete if we can show that, for any $L \in {\mathbb N}$,
\beq \label{but2pthmItem2} \Vert D_{L,1} \Vert_1=  O (2^{L/p} L^{a_p} ) \ \, \text{  and } \ \, \Vert D_{L,2} \Vert_1=  O (2^{L/p} L^{a_p} )  \, . 
\eeq
As we shall see below, this will be achieved by selecting  the sequence $(m(L))_{L \geq 0}$ as follows in the construction of the iid gaussian vectors $(Y_i)_{i \geq 1}$ as described previously: set 
$b_p = \frac{1}{p} $  if $p \in ]2, 3[$  and $b_3= -1/3 $ and 
\beq \label{defmlitem2}
m(L) = \Big [ \frac{2L}{p} - b_p \log_2 L \Big ] \, ,\ \text{so that}\ 
\frac 12 2^{2L/p}L^{-b_p}  \leq 2^{m(L)} \leq 2^{2L/p} L^{-b_p} \, .
\eeq 
The first part  of \eqref{but2pthmItem2} follows by using  \eqref{cons1keylemma} together with the above selection of $m(L)$. 
To show the second part of \eqref{but2pthmItem2}, we set $y_L = \kappa_m^{-1/2} ( \log 2 )^{1/2} 2^{m(L) /2} L^{1/2}$ where $\kappa_m = \lambda_m /c(m) $. Hence, using \eqref{decdl2bis}, write that 
\begin{align*}
\Vert D_{L,2} \Vert_1 & \leq y_L  + \int_{y_L }^{\infty} {\mathbb P} ( D_{L,2} \geq t) dt  \\
& \leq C \Big \{  y_L  + 2^L y_L^{1-p} +  2^{L-m(L)} \int_{y_L }^{\infty} \exp \Bigl( - \frac{\kappa_m t^2  }{      2^{m(L) }} \Bigr)   dt \Big \}  \\
& \leq C \Big \{ y_L  + 2^L y_L^{1-p} +  2^{L}  \exp \Bigl( - \frac{\kappa_m y_L^2  }{      2^{m(L) }} \Bigr)  \Big \}   \leq C \{ y_L  + 2^L y_L^{1-p}  \} \, .
\end{align*}
Taking into account the selection of $y_L$ and \eqref{defmlitem2}, Item 2 follows.

%%%%%%%%%%%%%%%%%%%%%%%%%%%%%%%%%

\subsubsection{Proof of Theorem \ref{ThmartW1cond}} 

Recall the notations  $A^{-1} = P  \Delta^{1/2}$ and $m_k = A d_k$. Hence, we have
 \[
 W_{1} (P_{M_n |{\mathcal F}_0}, P_{G_{n \Sigma}} ) \leq \Vert A^{-1} \Vert_2  W_{1} (P_{ \sum_{k=1}^n m_k  |{\mathcal F}_0 }, P_{G_{n  {\rm I}_m}} ) = \sqrt{\lambda_1} W_{1} (P_{ \sum_{k=1}^n m_k  |{\mathcal F}_0}, P_{G_{n J_m}} )  \text{ a.s.}
\]
Moreover, since $ \E (m_0 m_0^t)= {\rm J}_m$, we have $(m_0)_i =0$ a.s. for any $i=m+1, \ldots, d$. Hence, setting 
\[
d_k' = \big ( (m_k)_1, \ldots,  (m_k)_m\big )^t \,  \text{ and } \,  M_n'= \sum_{k=1}^n d'_k\, , 
\]
and noticing that $ \E (d'_0 (d'_0)^t)= {\rm I}_m$, we have
\[
W_{1} (P_{ \sum_{k=1}^n m_k  |{\mathcal F}_0}, P_{G_{n J_m}} ) = W_{1} (P_{ M'_n |{\mathcal F}_0}, P_{G_{n I_m}} ) \text{ a.s.}
\]
From Point 2 of Theorem 1 in \cite{DPR}, the following inequality  then  holds: 
\[
\Vert  W_{1} (P_{M_n |{\mathcal F}_0}, P_{G_{n \Sigma}} ) \Vert_1  \leq  \sqrt{\lambda_1}  \sup_{g \in {\rm Lip}(|\cdot|_{m}, {\mathcal F}_{0})} {\mathbb E}  (  g(M'_n )  ) 
- {\mathbb E} (g( T_n )) \, ,
\]
where $T_n = \sum_{i=1}^n N_i$ with $(N_i)_{i
\geq 1}$ a sequence of iid ${\mathbb R}^m$-valued centered gaussian random variables with ${\rm Var} (N_1) = {\rm I}_m$. To end the proof of the theorem, it suffices to use Lemma \ref{keylemma}.

%%%%%%%%%%%%%%%%%%%%%%%%%%%%%%%%%

\subsubsection{Proof of Lemma \ref{keylemma}} Note first that we can assume the sequence $(N_i)_{i
\geq 1}$ independent of $(d_i)_{i \in {\mathbb Z}}$. Let us now consider a $m$-dimensional standard Gaussian random vector $G$ independent of 
$(N_i)_{i
\geq 1}$ and $(d_i)_{i \in {\mathbb Z}}$. Note that 
\begin{align*} 
\sup_{g \in {\rm Lip}( |\cdot |_{m}, {\mathcal F}_{0})} {\mathbb E}  (  g(M'_n)  ) 
- {\mathbb E} (g( T_n)) & \leq \sup_{g \in {\rm Lip}( |\cdot |_{m}, {\mathcal F}_{0})} {\mathbb E}  (  g(M'_n +G )  ) 
- {\mathbb E} (g( T_n + G))  +  2  {\mathbb E} |G|_m \\
&  \leq  \sup_{g \in {\rm Lip}( |\cdot |_{m}, {\mathcal F}_{0})} {\mathbb E}  (  g(M'_n +G )  ) 
- {\mathbb E} (g( T_n + G)) + 2  \sqrt{m}  \, .
\end{align*}
The lemma is then reduced to prove that 
\beq \label{p1keylemma}
 \sup_{g \in {\rm Lip}( |\cdot |_{m}, {\mathcal F}_{0})} {\mathbb E}  (  g(M'_n +G )  ) 
- {\mathbb E} (g( T_n + G))  \leq   \left\{
  \begin{array}{rcr}
   C  n^{ (3-p) /2} & \text{ if $p \in ]2,3[$} \\
  C \log n  & \text{ if $p=3 $} \\
  \end{array}
\right.\, . 
\eeq
We now use the Lindeberg method to prove \eqref{p1keylemma}.  With this aim, we introduce the following notation:

\begin{nota}\label{not21} Let $\varphi_a$ be the density of a $m$-dimensional centered Gaussian random vector with covariance matrix $a^2 {\rm I}_m$ and let for $x  \in {\mathbb R}^m$,
$$g * \varphi_a (x,\omega) = \int  g ( x+y, \omega ) \varphi_a(y) dy\, . $$ 
For the sake of brevity, we shall write $g * \varphi_a (x)$ 
instead of  $g * \varphi_a (x,\omega)$ (the partial derivatives will be taken with respect to $x$) and set $g_{i,n} (x) =g * \varphi_{i,n}$  where $\varphi_{i,n}=\varphi_{\sqrt{n -i +1}}$. 
\end{nota}
Let
\[
\Delta_{i,n} (g)  =  g\big (M'_{i-1} + d'_{i}+ \sum_{j= i+1}^{n} N_j + G\big ) - g\big (M'_{i-1} + N_{i} +  \sum_{j= i+1}^{n} N_j + G\big )    \, .
\]
We have
\[
\E (\Delta_{i,n} (g)) = \E \big (g_{i,n}\big (M'_{i-1} + d'_{i} \big ) \big )  -\E \big ( g_{i,n}\big (M'_{i-1} + N_{i}  \big )  \big ) 
\]
Hence, noticing that $M_0 =0$,   the following decomposition is valid: 
\beq \label{lind1}
 {\mathbb E}  (  g(M'_n +G)  ) 
- {\mathbb E} (g( T_n  + G))  =  \sum_{i=1}^{n} \bkE \big ( \Delta_{i,n} (g) \big )  \, . 
\eeq

Below we shall also use the following notations.

\begin{nota} \label{deftensorprod} For two positive integers $m$ and $n$, let ${\mathcal M}_{m,n} ({\mathbb R})$ be the set of real matrices with $m$ lines and $n$ columns. The Kronecker product  
(or Tensor product) of $A =  [ a_{i,j} ] \in {\mathcal M}_{m,n} ({\mathbb R})$ and $B =  [ b_{i,j} ] \in {\mathcal M}_{p,q} ({\mathbb R})$ is denoted by $A \otimes B$ and is defined to be the block matrix 
$$
A \otimes B = \left ( \begin{array}{ccc}
a_{1,1}B & \cdots & a_{1,n}B \\ 
\vdots &  & \vdots \\ 
a_{m,1}B & \cdots & a_{m,n}B
\end{array} \right ) \in {\mathcal M}_{mp,nq} ({\mathbb R}) \, .
$$
For any positive integer $k$, the $k$-th Kronecker power $A^{\otimes k}$ is defined inductively by: $A^{\otimes 1}=A$ and $A^{\otimes k} = A \otimes A^{\otimes (k-1)}$. 

If $\nabla$ denotes the differentiation operator given by $\nabla  =   \big ( \frac{\partial }{\partial x_1}, \ldots, \frac{\partial }{\partial x_m} \big )^t$ acting on  the differentiable functions $f: {\mathbb R}^m \rightarrow {\mathbb R}$,   we define 
$$
\nabla \otimes \nabla =  \Big ( \frac{\partial }{\partial x_1} \circ \nabla, \ldots, \frac{\partial }{\partial x_m} \circ \nabla \Big )^t\, ,
$$
and  $\nabla^{\otimes k}$  by $\nabla^{\otimes 1}=\nabla$ and $\nabla^{\otimes k} = \nabla \otimes \nabla^{\otimes (k-1)}$. If $f : {\mathbb R}^m \rightarrow {\mathbb R}$ is $k$-times differentiable, for any $x \in {\mathbb R}^m$, let $
D^k f(x) = \nabla^{\otimes k} f(x) $, 
and for any vector $A$ of ${\mathbb R}^m$, we define $D^k f(x) \text{{\bf .}} A^{\otimes k}$ as the usual scalar product in ${\mathbb R}^{m^k}$ between $D^k f(x)$ and $ A^{\otimes k}$. 
\end{nota}
For any $i\in \{1, \ldots, n \}$, let  
$$
\Delta_{1,i,n} (g ) =  g_{i,n}\big (M'_{i-1} + d'_{i} \big )  -  g_{i,n} \big (M'_{i-1} \big ) 
- \frac{1}{2} D^2 g_{i,n} \big (M'_{i-1}\big ) \text{{\bf .}}   N_{i}^{\otimes 2} \, ,
$$
and 
$$
\Delta_{2,i,n} (g) = g_{i,n}\big (M'_{i-1} + N_{i} \big )  -  g_{i,n} \big (M'_{i-1} \big ) 
- \frac{1}{2} D^2 g_{i,n} \big (M'_{i-1}\big ) \text{{\bf .}}   N_{i}^{\otimes 2}  \, .
$$
With this notation,
\begin{equation} \label{lind2}
 \bkE \big ( \Delta_{i,n} (g) \big ) = \bkE \big ( \Delta_{1,i,n} (g) \big ) - \bkE \big (\Delta_{2,i,n} (g)  \big ) \, . 
\end{equation}
By the Taylor integral formula, noticing that $\bkE (  N_{i}^{\otimes 3}) = 0$, we get 
$$
\big | \bkE (\Delta_{2,i,n} (g) )\big | \leq \frac{1}{6}  \Big | \bkE  \int_0^1 D^4 g_{i,n} \big (M'_{i-1} + t  N_{i} \big ) \text{{\bf .}}   N_{i}^{\otimes 4} dt \Big | \, .
$$
But, according to Lemma 5.6 in \cite{DMR13}, for any $y \in {\mathbb R}^{m}$ and  any integer $k \geq 1$,  there exists a positive  constant $c_k$ depending only on $k$ such that 
\beq \label{Lemma56DMR}
\sup_{(i_1, \ldots, i_k) \in \{1, \ldots, m \}^m}\Big |  \frac{\partial^k  g * \varphi_{i,n}}{\prod_{j=1}^k\partial x_{i_j}} (y) \Big | \leq c_k (n -i +1 )^{(1-k)/2 } \, .
\eeq
Therefore
\begin{equation*} 
\big | \bkE (\Delta_{2,i,n} (g) )\big | \leq 2^{-1} m^4 c_4 (n -i +1 )^{-3/2 } \, .
\end{equation*}
Therefore
\beq \label{lind3}
\sum_{i=1}^{n} \big | \bkE (\Delta_{2,i,n} (g) )\big | 
\leq 2^{-1} m^4 c_4   \sum_{i=1}^{n }  (n -i +1 )^{-3/2 } \leq \kappa_1 m^4 \, .
\eeq
Let now 
$$
R_{1,i,n} (g ) = g_{i,n}  \big (M'_{i-1} + d'_{i} \big ) -g_{i,n}  \big (M'_{i-1} \big ) -  D g_{i,n}  \big (M'_{i-1}\big ) \text{{\bf .}}  d'_{i} 
- \frac{1}{2} D^2 g_{i,n}  \big (M'_{i-1}\big ) \text{{\bf .}}  d_{i}^{ \prime \otimes 2} \, ,
$$
and 
$$
R_{2,i,n} (g) =  D g_{i,n}  \big (M'_{i-1}\big ) \text{{\bf .}}  d'_{i} +  \frac{1}{2} D^2 g_{i,n}  \big (M'_{i-1}\big ) \text{{\bf .}}  d_{i}^{\prime \otimes 2} - \frac{1}{2} D^2 g_{i,n}  \big (M'_{i-1}\big ) \text{{\bf .}}
 \bkE ( N_{i}^{\otimes 2} )\, .
$$
With this notation,
\begin{equation} \label{lind4}
\bkE (\Delta_{1,i,n} (g) ) = \bkE (R_{1,i,n} (g) ) + \bkE (R_{2,i,n} (g) ) \, .
\end{equation}
By  the Taylor integral formula at order two,
$$
\big |\bkE (R_{1,i,n} (g) )\big | \leq  \Big |  \bkE \int_0^1 \frac{ (1-t) }{2} \Big (  D^2 g_{i,n} \big (M'_{i-1} + t d'_i \big ) -   D^2 g_{i,n} \big (M'_{i-1}  \big )  \Big )  \text{{\bf .}}  d_{i}^{ \prime \otimes 2} dt  \Big | \, .
$$
But, by taking into account \eqref{Lemma56DMR}, we infer that, for any $t \in [0,1]$, 
\begin{multline*}
\big | D^2 g_{i,n} \big (M'_{i-1} + t d'_i \big ) -   D^2 g_{i,n} \big (M'_{i-1}  \big )  \Big )  \text{{\bf .}}  d_{i}^{\prime \otimes 2} \big | \\
 \leq 
m (c_2 \vee c_3) (n -i +1 )^{ -1/2 } \Big ( \min \big ( 2  | d'_1 |_m^2    ,  t  m^{1/2}  (n -i +1 )^{-1/2} | d'_1 |_m^3 \big )  \Big )  \\
 \leq 
2 m^{p/2} (c_2 \vee c_3) (n -i +1 )^{ -1/2 } (n -i +1 )^{-(p-2)/2}  |  d'_1 |_m^p \, ,
\end{multline*}
since $p \in ]2,3]$. Therefore, since $  |  d'_1 |_m = | A d_1 |_d \leq \Vert A \Vert_2 | d_1 |_d = \sqrt{\rho(P \Gamma^2 P^t) } | d_1 |_d = \lambda_m^{-1/2}  | d_1 |_d $, it follows that 
\[
\big |\bkE (R_{1,i,n} (g) )\big | \leq   2^{-1} \lambda_m^{-p/2} m^{p/2}   (c_2 \vee c_3)   (n -i +1 )^{-(p-1)/2}  \bkE  (  |  d_1 |_d^p )  \, .
\]
Hence there exists a positive constant $\kappa_2$ such that, for any positive integer $n$, 
\beq \label{lind5}
\sum_{i=1}^{n} \big |\bkE (R_{1,i,n} (g) )\big | 
\leq \kappa_2  \lambda_m^{-p/2} m^{p/2}  \bkE  (  |  d_1 |_d^p ) n^{(3-p)/2}   ( 1 + {\bf 1}_{p=3} \log n  ) \, .
\eeq
Starting from the decomposition \eqref{lind2} and taking into account \eqref{lind3},  \eqref{lind4} and  \eqref{lind5}, it follows that for any $g \in {\rm Lip}(|\cdot|_{m}, {\mathcal F}_{0})$, 
\beq \label{lind6}
\Big |  \sum_{i=1}^{n} \bkE \big ( \Delta_{i,n} (g) \big )  \Big | \leq \Big |  \sum_{i=1}^{n}  (\bkE (R_{2,i,n} (g) )   \Big |  + C n^{(3-p)/2}   ( 1 + {\bf 1}_{p=3} \log n  ) \, .
\eeq

We handle now the term $ \sum_{i=1}^{n}\bkE (R_{2,i,n} (g) ) $.  With this aim, let us first write $n$ in basis $2$. Let $r$ be the unique non-negative integer such that $2^r \leq n < 2^{r+1}$. Then, writing $n$ in basis $2$, we have
\[
n=\sum_{k=0}^r b_k(n) 2^k  \ \text{ where $b_r(n)=1$ and $b_k(n)\in \{0, 1 \}$ for $k =0, \ldots, r-1$} \, .
\]
Let also
\[
n_k = \sum_{j=0}^k b_j(n) 2^j \ \text{ for $k =0, \ldots, r-1$ and } n_{-1}=0 \, .
\]
It follows that 
\beq \label{lind1*1}
  \sum_{i=1}^{n}\bkE (R_{2,i,n} (g) )  =  \sum_{k=0}^r b_k(n) \sum_{i=n_{k-1} +1}^{n_k} \bkE (R_{2,i,n} (g) ) =  \sum_{k=0}^r b_k(n) \sum_{i=1}^{2^k} \bkE (A_{i,n} (g) )  \, , 
\eeq
where $A_{i,n} (g)  = R_{2,i+n_{k-1},n} (g)$  (note that the last inequality above holds because, when $b_k(n) =1$, then necessarily $n_k = 2^k + n_{k-1}$). 
Let $i \in \{1, \ldots, 2^k \}$. Notice now that, since $(d'_{i})_{i \in {\mathbb Z}} $ is 
a martingale differences sequence with respect to $({\mathcal F}_{i})_{i \in {\mathbb Z}} $ and such that 
$ \bkE ( d_{i}^{\prime \otimes 2} ) =  \bkE ( N_{i}^{\otimes 2} )$, we have, by setting $i(k) = i+n_{k-1}$, 
\[
 \bkE (A_{i,n} (g) )=  \frac{1}{2}  \bkE \big ( D^2 g_{i(k),n}  \big (M'_{i(k)-1}\big ) \text{{\bf .}}   \big ( d_{i(k)}^{ \prime \otimes 2}  - \bkE ( d_{i(k)}^{ \prime \otimes 2} )   \big ) \, .
\]
To continue the computations, as in the proof of Proposition 5.1 in \cite{DMR09}, we introduce again a dyadic
scheme. With this aim, we introduce the following notations. 
\begin{nota} Set  $i_0 = i-1$ and write $i_0$ in basis $2$ (recall that $i \in \{1, \ldots, 2^k \}$): $ i_0 = \sum_{i=0}^{k}
a_i 2^i $ with  $a_i=0$ or $a_i=1$ (note that $a_{k}  = 0$). Set $ i_j
=   \sum_{i=j}^{k} a_i 2^i , $ so that $i_{k} = 0$, and set $i_j(k) = i_j + n_{k-1}$. Let $I_{j,\ell} =  ]\ell
2^j , (\ell+1) 2^j] \cap {\mathbb N}$ (note that $I_{k,1} = ]2^k,
2^{k+1}]$), $U_j^{(\ell)} = \sum_{i \in I_{j,\ell}  }  d'_{i }$, $U_{j,k}^{(\ell)} = \sum_{i \in I_{j,\ell}  }  d'_{i +k}$,  $\tilde
U_j^{(\ell)} = \sum_{i \in I_{j,\ell}  }  N_{i }$ and $\tilde
U_{j,k}^{(\ell)} = \sum_{i \in I_{j,\ell}  }  N_{i +k}$. For the sake of brevity,
let $U_{j}^{(0)}=U_j$ and $\tilde U_{j}^{(0)}= \tilde U_j$. Set also $Z_j^{(\ell)} = \E_{\ell 2^j } (( U_j^{(\ell)} )^{\otimes 2}) -  \E_{\ell 2^j } (( U_j^{(\ell)} )^{\otimes 2})$. 
\end{nota}

Since  $i_{k} =0$, the following elementary identity is valid
\begin{multline} \label{lind7}
 \bkE (A_{i,n} (g) ) \\
 =  \frac{1}{2}  \sum_{j=0}^{k-1} \bkE \Bigl(  \big ( D^2 g_{i_j (k)+1,n}  \big (M'_{i_j (k)}\big )     -   D^2 g_{i_{j+1}(k)+1,n}  \big (M'_{i_{j+1} (k)} \big ) \big )  \text{{\bf .}}  
\big ( d_{i(k)}^{\prime \otimes 2}  - \bkE ( d_{i(k)}^{ \prime \otimes 2} )   \big ) \Bigr) \\
+  \frac{1}{2}  \bkE \Bigl(    D^2 g_{n_{k-1} + 1,n}  \big ( M'_{n_{k-1}} \big )   \text{{\bf .}}  \big ( d_{i(k)}^{ \prime \otimes 2}  - \bkE ( d_{i(k)}^{ \prime \otimes 2} )   \big ) \Bigr)  \\
=: \bkE (B_{i,n} (g) )  + \frac{1}{2}  \bkE \Bigl(    D^2 g_{n_{k-1} + 1,n}  \big ( M'_{n_{k-1}} \big )   \text{{\bf .}}  \big ( d_{i(k)}^{\prime \otimes 2}  - \bkE ( d_{i(k)}^{ \prime \otimes 2} )   \big )   \Bigr) \, . 
\end{multline}
Since $ D^2 g_{n_{k-1} + 1,n}  \big ( M'_{n_{k-1}} \big )  $ is a random vector which is ${\mathcal F}_{n_{k-1}}$ measurable, we have
\[
\bkE \Bigl(     D^2 g _{n_{k-1} + 1,n}  \big ( M'_{n_{k-1}} \big )  \text{{\bf .}}    \big ( d_{i(k)}^{\prime \otimes 2}  - \bkE ( d_{i(k)}^{\prime \otimes 2} )   \big )   \Bigr)   = 
\bkE \Bigl(    D^2 g_{n_{k-1} + 1,n}  \big ( M'_{n_{k-1}} \big )    \text{{\bf .}}  \bkE_{n_{k-1}}  \big ( d_{i(k)}^{\prime \otimes 2}  - \bkE ( d_{i(k)}^{\prime \otimes 2} )   \big )  \Bigr) \, .
\]
But, by the martingale property,
\[
\sum_{i=1}^{2^k}  \bkE_{n_{k-1}}   \big ( d_{i(k)}^{\prime \otimes 2}  - \bkE ( d_{i(k)}^{\prime \otimes 2} )   \big ) =  \bkE_{n_{k-1}}  \Big ( \Big ( \sum_{i=1}^{2^k}  d'_{i(k)}  \Big)^{\otimes 2}  - \bkE  \Big ( \sum_{i=1}^{2^k}   d'_{i(k)} \Big )^{\otimes 2}   \Big ) 
\, .
\]
Taking into account \eqref{Lemma56DMR} and stationarity, we get
\begin{multline*}
\Big | \sum_{i=1}^{2^k} \bkE \Bigl(     D^2 g _{n_{k-1} + 1,n}  \big ( M'_{n_{k-1}} \big )  \text{{\bf .}}    \big ( d_{i(k)}^{\prime \otimes 2}  - \bkE ( d_{i(k)}^{\prime \otimes 2} )   \big )   \Bigr) \Big | \\
 \leq  c_1   (n-n_{k-1})^{-1/2}  \sum_{a=1}^m \sum_{b =1}^m  \bkE \Big (   \Big |  \bkE_{0} \Big ( \Big ( \sum_{i=1}^{2^k} d'_{i} \Big )_a \Big ( \sum_{i=1}^{2^k} d'_{i} \Big )_{b} \Big ) - 
\bkE\Big ( \Big ( \sum_{i=1}^{2^k} d'_{i} \Big )_a \Big ( \sum_{i=1}^{2^k} d'_{i} \Big )_{b} \Big )   \Big |  \Big ) \, .
\end{multline*}
Since $(d'_i)_{a} = (Ad_i)_a $ for $1 \leq a \leq m$ and $ (Ad_i)_a = 0$ for $m+1 \leq a \leq d$, 
\begin{multline*}
  \sum_{a=1}^m \sum_{b =1}^m  \bkE \Big (   \Big |  \bkE_{0} \Big ( \Big ( \sum_{i=1}^{2^k} d'_{i} \Big )_a \Big ( \sum_{i=1}^{2^k} d'_{i} \Big )_{b} \Big ) - 
\bkE\Big ( \Big ( \sum_{i=1}^{2^k} d'_{i} \Big )_a \Big ( \sum_{i=1}^{n} d'_{i} \Big )_{b} \Big )   \Big |  \Big )  \\
\leq d  \,  \bkE \Big |   A^{\otimes 2} \big (    \bkE_{0} \big (  M_{2^k}^{\otimes 2} \big ) -  \bkE \big (  M_{2^k}^{\otimes 2} \big )   \big )  \Big |_d \leq  d   \Vert  A^{\otimes 2} \Vert_2  \bkE \Big |     \bkE_{0} \big (  M_{2^k}^{\otimes 2} \big ) -  \bkE \big (  M_{2^k}^{\otimes 2} \big )    \Big |_d  \, .
\end{multline*}
But $
\Vert  A^{\otimes 2} \Vert_2  \leq \Vert  A\Vert^2_2 \leq \lambda_m^{-1} $. So, from the above considerations, 
\begin{multline*}
\Big | \sum_{i=1}^{2^k} \bkE \Bigl(     D^2 g_{n_{k-1} + 1,n}  \big ( M'_{n_{k-1}} \big )  \text{{\bf .}}    \big ( d_{i(k)}^{\prime \otimes 2}  - \bkE ( d_{i(k)}^{\prime \otimes 2} )   \big )   \Bigr) \Big | \\
 \leq  c_1  d \lambda_m^{-1}   (n-n_{k-1})^{-1/2}  \bkE \Big |     \bkE_{0} \big (  M_{2^k}^{\otimes 2} \big ) -  \bkE \big (  M_{2^k}^{\otimes 2} \big )    \Big |_d  \, .
\end{multline*}
Since $n \geq n_{k-1} +2^k $ if $b_k(n) =1$, we get 
\begin{multline*} 
\sum_{k=0}^r b_k(n)\Big | \sum_{i=1}^{2^k} \bkE \Bigl(     D^2 g_{n_{k-1} + 1,n}  \big ( M'_{n_{k-1}} \big )  \text{{\bf .}}    \big ( d_{i(k)}^{\prime \otimes 2}  - \bkE ( d_{i(k)}^{\prime \otimes 2} )   \big )   \Bigr) \Big | \\
\leq C  \sum_{k=0}^r b_k(n) 2^{-k/2}  \bkE \Big |     \bkE_{0} \big (  M_{2^k}^{\otimes 2} \big ) -  \bkE \big (  M_{2^k}^{\otimes 2} \big )    \Big |_d   \\
\leq C 2^{ (3-p) r /2}   \sum_{k=0}^r   2^{ - \big ( 2-\frac{p}{2} \big ) k}  \bkE \Big |     \bkE_{0} \big (  M_{2^k}^{\otimes 2} \big ) -  \bkE \big (  M_{2^k}^{\otimes 2} \big )    \Big |_d   \, .\end{multline*}
But, by the subadditivity of the sequence $ \big (  \bkE \big |     \bkE_{0} \big (  M_{k}^{\otimes 2} \big ) -  \bkE \big (  M_{k}^{\otimes 2} \big )    \big |_d \big )_{k \geq 1}$, condition \eqref{C1} implies that 
\[
 \sum_{k\geq 0}   2^{ - \big ( 2-\frac{p}{2} \big ) k}  \bkE \Big |     \bkE_{0} \big (  M_{2^k}^{\otimes 2} \big ) -  \bkE \big (  M_{2^k}^{\otimes 2} \big )    \Big |_d  < \infty \, .
\]
(See for instance Remark 2.6 in \cite{DMR09}). Hence, 
\beq \label{lind8}
\sum_{k=0}^r b_k(n)\Big | \sum_{i=1}^{2^k} \bkE \Bigl(     D^2 g_{n_{k-1} + 1,n}  \big ( M'_{n_{k-1}} \big )  \text{{\bf .}}    \big ( d_{i(k)}^{\prime \otimes 2}  - \bkE ( d_{i(k)}^{\prime \otimes 2} )   \big )   \Bigr) \Big | \leq C n^{(3-p)/2} \, .\eeq
Starting from \eqref{lind6} and taking into account \eqref{lind1*1}, \eqref{lind7} and \eqref{lind8}, the proof of the lemma will be complete if we can prove that 
for any $g \in {\rm Lip}(|\cdot |_m, {\mathcal F}_{0})$, 
\beq \label{lind9}
 \Big |  \sum_{k=0}^r b_k(n) \sum_{i=1}^{2^k} \bkE (B_{i,n} (g) )  \Big |  \leq  C n^{(3-p)/2}( 1 + {\bf 1}_{p=3} \log n  ) \, .
\eeq
Let $i \in \{1, \ldots, 2^k \}$ and  note that $i_j \neq i_{j+1}$ only if $a_j = 1$, then in this case $i_j=\ell
2^j$ with $\ell$ odd. It follows that 
\begin{multline*} 
2 \sum_{i=1}^{2^k} \bkE (B_{i,n} (g) ) 
 =   \sum_{j=0}^{k-1} \sum_{\ell
\in I_{k-j ,0} \atop \ell \text { odd }}  \bkE \Bigl(  \big ( D^2 g_{\ell
2^j + n_{k-1}+1,n}  \big (M'_{\ell
2^j+n_{k-1}}\big )    \\  -   D^2 g_{(\ell -1)
2^j + n_{k-1}+1,n}  \big (M'_{(\ell-1)
2^j+n_{k-1}}\big ) \big )  
 \text{{\bf .}}  
\sum_{\{i : i_j = \ell 2^j\}} \big ( d_{i(k)}^{\prime \otimes 2}  - \bkE ( d_{i(k)}^{\prime \otimes 2} )   \big ) \Bigr)    \, . 
\end{multline*}
Note
that $\{i : i_j = \ell 2^j\}=I_{j,\ell}$. Now by the martingale property,
$$ {\mathbb E}_{\ell 2^j +n_{k-1}}\Big(\sum_{ i\in I_{j,\ell}}
 \big ( d_{i(k)}^{\prime \otimes 2}  - \bkE ( d_{i(k)}^{\prime \otimes 2} )   \big ) \Big)=  {\mathbb E}_{\ell 2^j+n_{k-1}} ( ( U_{j,n_{k-1}}^{(\ell)} )^{\otimes 2}  ) -
\bkE ( ( U_{j,n_{k-1}}^{(\ell)} )^{\otimes 2}  ):=Z_{j,n_{k-1}}^{(\ell)}\, .
$$
Consequently
\begin{align*}
2 \sum_{i=1}^{2^k} \bkE (B_{i,n} (g) )  & = \sum_{j=0}^{k-1} \sum_{\ell
\in I_{k-j ,0} \atop \ell \text { odd }}  \bkE \Bigl(  \big ( D^2 g_{\ell
2^j + n_{k-1}+1,n}  \big (M'_{\ell
2^j+n_{k-1}}\big )     \\ & \quad \quad -   D^2 g_{(\ell -1)
2^j + n_{k-1}+1,n}  \big (M'_{(\ell-1)
2^j+n_{k-1}}\big ) \big )  \text{{\bf .}} 
Z_{j,n_{k-1}}^{(\ell)} \Bigr) \\
& = \sum_{j=0}^{k-1} \sum_{\ell
\in I_{k-j ,0} \atop \ell \text { odd }}  \bkE \Bigl(  \Big ( D^2 g _{\ell
2^j + n_{k-1}+1,n}  \big (M'_{\ell
2^j+n_{k-1}}\big )    \\ & \quad \quad -     D^2 g_{\ell 
2^j + n_{k-1}+1,n}  \Big (M'_{(\ell-1)
2^j+n_{k-1}} +  \sum_{v= (\ell-1) 2^j +1}^{\ell 2^j} N_{v + n_{k-1}} \Big ) \Big )  \text{{\bf .}} 
Z_{j,n_{k-1}}^{(\ell)} \Bigr) 
\, ,
\end{align*}
since $(d_i)_{i \in {\mathbb N}}$ and  $(N_i)_{i \in {\mathbb N}}$ 
are independent. Note that $Z_{j, n_{k-1}}^{(\ell)}$ is a $m^2$-dimensional random vector. Let then introduce the following notation: $(Z_{j,n_{k-1}}^{(\ell)})_{a,b}$ is the $(a-1)m +b$-th coordinate of the vector 
$Z_{j,n_{k-1}}^{(\ell)}$. By using this notation and \eqref{Lemma56DMR}, it follows that 
\begin{multline*}
2 \Big | \sum_{i=1}^{2^k} \bkE (B_{i,n} (g) )  \Big |  \leq  \max (2c_2,c_3)  \times  \sum_{j=0}^{k-1} \sum_{\ell
\in I_{k-j ,0} \atop \ell \text { odd }} 
 (n - \ell 2^j -n_{k-1})^{-1/2}   \\
 \times \sum_{a,b=1}^m  
\bkE \Bigl(  
| (Z_{j,n_{k-1}}^{(\ell)} )_{a,b}  | \min  \big ( 1 ,    (n - \ell 2^j -n_{k-1})^{-1/2}  \sum_{c=1}^m \big | (   U_{j,n_{k-1}}^{(\ell-1)} - \ti U_{j,n_{k-1}}^{(\ell-1)} )_c  \big |  \big )  \Bigr)
\, .
\end{multline*}
But, since $p \in ]2,3]$, 
\begin{multline*}
\bkE \Bigl(  
| (Z_{j,n_{k-1}}^{(\ell)} )_{a,b}  | \min  \big ( 1 ,    (n - \ell 2^j -n_{k-1})^{-1/2}  \sum_{c=1}^m \big | (   U_{j,n_{k-1}}^{(\ell-1)} - \ti U_{j,n_{k-1}}^{(\ell-1)} )_c  \big |  \big )  \Bigr)
 \\
\leq  (n - \ell 2^j -n_{k-1})^{-(p-2)/2}   \bkE \Bigl(  
|  (Z_{j,n_{k-1}}^{(\ell)} )_{a,b}  | \Big (  \sum_{c=1}^m \big | (   U_{j,n_{k-1}}^{(\ell-1)} - \ti U_{j,n_{k-1}}^{(\ell-1)} )_c  \big |  \Big )^{p-2}  \Bigr) \\
\leq  (n - \ell 2^j -n_{k-1} )^{-(p-2)/2}   \sum_{c=1}^m  \bkE \Bigl(  
|  (Z_{j,n_{k-1}}^{(\ell)} )_{a,b}   \big | (   U_{j,n_{k-1}}^{(\ell-1)} - \ti U_{j,n_{k-1}}^{(\ell-1)} )_c  \big |^{p-2}  \Bigr)
\, ,
\end{multline*}
which together with stationarity implies that 
\begin{multline*}
2 b_k(n)  \Big | \sum_{i=1}^{2^k} \bkE (B_{i,n} (g) )  \Big |  \leq \\  \leq  \max (2c_2,c_3)  \sum_{j=0}^{k-1} \sum_{\ell =1}
^{2^{k-j} -1}
  b_k(n) (n - \ell 2^j -n_{k-1})^{-(p-1)/2} 
 \sum_{a,b,c=1}^m  
\bkE \Bigl(  
| (Z_{j}^{(1)} )_{a,b} | \big | (   U_{j} - \ti U_{j} )_c  \big |^{p-2}  \Bigr) 
\, .
\end{multline*}
Let
\[
Q_j =  \sum_{a,b,c=1}^m  
\bkE \Bigl(  
| (Z_{j}^{(1)} )_{a,b} | \big | (   U_{j} - \ti U_{j} )_c  \big |^{p-2}  \Bigr) \, .
\]
Since  $n \geq n_{k-1} +2^k $ if $b_k(n) =1$, we get, in case where $p \in ]2, 3[$, 
\[
b_k(n)  \Big | \sum_{i=1}^{2^k} \bkE (B_{i,n} (g) )  \Big |  
 \leq   C 
   \sum_{j=0}^{k-1}    2^{-j(p-1)/2}  2^{ (3-p) (k -j) /2}   Q_j \, .
\]
Consider now the case $p=3$. Note that, since $b_k(n) =1$ implies that $n_k = n_{k-1} +2^k $, we have 
\begin{multline*}
  b_k(n)  \sum_{\ell =1}
^{2^{k-j} -1}
  (n - \ell 2^j -n_{k-1})^{-1} =  2^{-j}b_k(n) \sum_{\ell = 2^{-j} (n-n_k) +1}
^{2^{-j} (n-n_{k-1}) -1}  \ell^{-1} \\  \leq 2^{-j +1}b_k(n) \sum_{\ell = 2^{-j} (n-n_k) +1}
^{2^{-j} (n-n_{k-1}) -1}  (\ell+1)^{-1} \leq  2^{-j +1}b_k(n)  \Big \{ \log \big ( \frac{n-n_{k-1}}{2^j}  \big ) - \log \big ( \frac{n-n_{k}}{2^j} +1 \big )   \Big \}  \, .
 \end{multline*}
 So, overall, 
\begin{multline*}
b_k(n)  \Big | \sum_{i=1}^{2^k} \bkE (B_{i,n} (g) )  \Big |  
 \leq   C  \left\{ \begin{array}{lll}
   \sum_{j=0}^{k-1}    2^{-j(p-1)/2}  2^{ (3-p) (k -j) /2}   Q_j  & \text{ if $p \in ]2,3[$} \\
     b_k(n)  \sum_{j=0}^{k-1}    2^{-j}   \Big \{ \log (n-n_{k-1} +1 ) - \log  (n-n_{k}+1  )   \Big \} Q_j  & \text{ if $p=3 $} \\
  \end{array}
\right.\, . 
 \end{multline*}
According to inequality (5.32) in \cite{DMR09}, we infer that under conditions \eqref{C1} and \eqref{C2}, 
\[ \bkE \Bigl(  
| (Z_{j}^{(1)} )_{a,b} | \big | (   U_{j} - \ti U_{j} )_c  \big |^{p-2}  \Bigr)   \leq  C 2^{ j ( 1- 2/p) }\Vert(Z_{j}^{(0)} )_{a,b} \Vert_{p/2} + C 2^{j(p/2-1)}
\Vert (Z_{j}^{(0)} )_{a,b} \Vert_{1, \Phi, p} \, .
\]
Hence, if $p \in ]2, 3[$, 
\begin{multline} \label{inefinallma}
 \sum_{k=0}^r b_k(n)  \Big | \sum_{i=1}^{2^k} \bkE (B_{i,n} (g) )  \Big | \\ \leq  C   \sum_{k=0}^r  2^{ (3-p) k /2}   \times 
 \sum_{a,b=1}^m   \Big \{     \sum_{j=0}^{k-1}    2^{-2j/p}  \Vert(Z_{j}^{(0)} )_{a,b} \Vert_{p/2}  +   \sum_{j=0}^{k-1} 
2^{j(p/2-2)} \Vert (Z_{j}^{(0)} )_{a,b} \Vert_{1, \Phi, p}   \Big \} 
\, ,
\end{multline}
and if $p=3$, 
\begin{multline} \label{inefinallmabis}
 \sum_{k=0}^r b_k(n)  \Big | \sum_{i=1}^{2^k} \bkE (B_{i,n} (g) )  \Big |   \leq  C   \sum_{k=0}^r   b_k(n)   \Big \{ \log (n-n_{k-1}  +1) - \log
 (n-n_{k}  +1 )   \Big \} \\
\times  \sum_{a,b=1}^m   \Big \{     \sum_{j=0}^{k-1}    2^{-2j/3}  \Vert(Z_{j}^{(0)} )_{a,b} \Vert_{3/2}  +   \sum_{j=0}^{k-1} 
2^{-j/2} \Vert (Z_{j}^{(0)} )_{a,b} \Vert_{1, \Phi, 3}   \Big \} 
\, .
\end{multline}
But, denoting by $ \alpha_{2^j , d} = \max_{k,\ell \in \{1, \ldots,  d \}}  \Vert \E_0 ( (  M_{2^j})_k  ( ( M_{2^j})_{\ell}   -  \E ( (M_{2^j})_k  ( ( M_{2^j})_{\ell} \Vert_{p/2} $, we get 
\begin{align*}
 \sum_{a,b=1}^m   \Vert(Z_{j}^{(0)} )_{a,b} \Vert_{p/2} & =  \sum_{a,b=1}^d  \Vert \E_0 ( ( A M_{2^j})_a  ( ( AM_{2^j})_{b}   -  \E ( (AM_{2^j})_a  ( ( AM_{2^j})_{b} \Vert_{p/2} \\
&  \leq   \sum_{a,b=1}^d \sum_{k , \ell =1}^d  \big | (A)_{a,k} (A)_{b,\ell} \big |\Vert \E_0 ( (  M_{2^j})_k  ( ( M_{2^j})_{\ell}   -  \E ( (M_{2^j})_k  ( ( M_{2^j})_{\ell} \Vert_{p/2}  \\
 &  \leq   d^2 \sum_{a,k=1}^d   (A)^2_{a,k}  \alpha_{2^j , d} := d^2 \Vert A \Vert_{HS}^2 \alpha_{2^j , d} \, .
\end{align*}
Now, since $P$ is orthogonal
\[
\Vert A \Vert^2_{HS} = \Vert \Gamma \Vert^2_{HS} =  \sum_{i=1}^{m-1}\lambda^{-1}_i+ (d-m+1) \lambda^{-1}_m \leq d \lambda^{-1}_m \, .
\]
So,
\[
\sum_{a,b=1}^m   \Vert(Z_{j}^{(0)} )_{a,b} \Vert_{p/2}
 \leq  d^3 \lambda^{-1}_m \alpha_{2^j , d} \, . 
 \]
Similarly, setting $ \beta_{2^j , d} = \max_{k,\ell \in \{1, \ldots,  d \}}  \Vert \E_0 ( (  M_{2^j})_k  ( ( M_{2^j})_{\ell}   -  \E ( (M_{2^j})_k  ( ( M_{2^j})_{\ell} \Vert_{1, \Phi, p} $, we get 
 \[
\sum_{a,b=1}^m   \Vert(Z_{j}^{(0)} )_{a,b} \Vert_{1, \Phi, p} 
  \leq d^3 \lambda^{-1}_m \beta_{2^j , d} \, . 
 \]
 Note now that due to the martingale property conditions \eqref{C1} and \eqref{C2} are respectively equivalent to 
\[
    \sum_{j \geq 0} 
2^{j(p/2-2)} \beta_{2^j , d}   < \infty \, \text{ and } \,  \sum_{j \geq 0}    2^{-2j/p}  \alpha_{2^j , d}  < \infty \, .
\]
(See Remark 2.6 in \cite{DMR09}). Moreover 
\[
\sum_{k=0}^r   b_k(n)   \Big \{ \log (n-n_{k-1}  +1) - \log
 (n-n_{k}  +1 )   \Big \} = \log (n+1) 
\]
and if $p \in ]2, 3[$,
\[
 \sum_{k=0}^r  2^{ (3-p) k /2}  \leq C n^{ (3-p) /2} 
\]
Starting from \eqref{inefinallma} and \eqref{inefinallmabis} and considering the above computations, the upper bound  \eqref{lind9} follows. This ends the proof of the  lemma. $\diamond$

\subsubsection{Proof of Remark \ref{commentreversed}} \label{proofcomment}

To simplify,  assume that $(d_i)_{i \in {\mathbb Z}}$ is a vector-valued stationary sequence of reversed martingale differences such that $\E (d_0d_0^t) = {\rm I}$ (otherwise we transform the r.v.'s to go back to this case as  done in Section \ref{sectpreli}). Let us then prove that if  \eqref{C1rev} and \eqref{C2rev} are satisfied for some $p \in ]2,3]$, the conclusions of Theorem \ref{Thmart} hold for $M_n$. With this aim, we need to construct the approximating Brownian motion. As in the proof of Theorem \ref{Thmart}, for $L \in {\mathbb N}$, let $m(L) \in {\mathbb N}$ be such that $m(L)\leq L$, and define 
\begin{equation*} 
I_{k,L} = ]2^L + (k-1)2^{m(L)} ,  2^L + k 2^{m(L)}] \cap {\mathbb N} \ \text{and}\ U_{k,L} = \sum_{i\in I_{k,L}} d_i \, ,
\, k \in \{1, \cdots, 2^{L-m(L)} \} \, .
\end{equation*}
Let $P_{U_{k,L} | {\mathcal G}_{2^L + k2^{m(L) +1}
}}$ be the conditional law of $U_{k,L}$ given $\mathcal{G}_{2^L +
k2^{m(L)} +1 }$ and ${\mathcal N}_{2^{m(L)}}$ be the ${\mathcal N} ( 0, 2^{m(L)} {\rm I})$-law. Assume that the probability space is assumed to be large enough to contain a sequence $(\delta_i)_{i \in {\mathbb Z}}$ of iid random variables uniformly distributed on $[0,1]$,  independent of the sequence $(d_i)_{ i \in {\mathbb Z}}$. As in the proof of Theorem \ref{Thmart},  we infer that 
there exists a ${\mathbb R}^d$-valued random variable $V_{k,L}$ with law ${\mathcal N}_{2^{m(L)}}$, measurable with respect to $\sigma(\delta_{2^L +k 2^{m(L)}}) \vee \sigma(U_{k ,L}) \vee {\mathcal G}_{2^{L}+ k2^{m(L)} +1} $, independent of ${\mathcal G}_{2^L + k2^{m(L)} +1}$ and such that 
\[
{\mathbb E} \big (  \big | U_{k,L} - V_{k,L} \big |_d \big ) = \sup_{g \in {\rm Lip}(|\cdot|_{d}, {\mathcal G}_{2^L + 2^{m(L)} +1})} {\mathbb E}  (  g(U_{1,L} )  ) 
- {\mathbb E} (g( V_{1,L}  )) \, .
\]
As, in the proof of Theorem \ref{Thmart}, by induction on $k$, we have then constructed Gaussian random variables $(V_{k,L})_{L \in {\mathbb N}, k=1, \ldots, 2^{L-m(L)}}$ that are mutually independent, and using the  Skorohod lemma, we can construct a sequence  $(Z_i)_{i \geq 1}$ of iid standard Gaussian random vectors in ${\mathbb R}^d$ such that, for any $L \in {\mathbb N}$ and any $ k \in \{1, \cdots, 2^{L-m(L)} \}$,
$ V_{k,L} =  \sum_{i\in I_{k,L}}  Z_i $  a.s.  The proof can be then completed if one can prove that, in  the case of reversed martingale differences, Lemma \ref{keylemma} still holds. In this case, it reads as follows: 
 \begin{Lemma} \label{keylemmarev}  Let $T_n = \sum_{i=1}^n N_i$ where $(N_i)_{i
\geq 1}$ be a sequence of iid ${\mathbb R}^d$-valued centered gaussian random variables with ${\rm Var} (N_1) = {\rm I}$. Let $p \in ]2,3]$ and assume that \eqref{C1rev} and \eqref{C2rev} are satisfied. Then, there exists a positive constant $C$ such that, for any positive integer $n$,  
\[
\sup_{g \in {\rm Lip}( |\cdot |_{d}, {\mathcal G}_{n+1})} {\mathbb E}  (  g(M_n)  ) 
- {\mathbb E} (g( T_n))  \leq   \left\{
  \begin{aligned}
  C  n^{ (3-p) /2} & \text{ if $p \in ]2,3[$} \\
  C \log n  & \text{ if $p=3 $.} \\
  \end{aligned}
\right.\, , 
\]
\end{Lemma} 
The proof of this lemma can be easily deduced from the one of Lemma \ref{keylemma} by writing $M_n = \sum_{j=1}^n d_{j,n}$ where $d_{j,n} = d_{n-j+1}$ for any  $j  \in {\mathbb Z}$ and by noticing that, for any integer $n$, 
$(d_{j,n})_{j \in {\mathbb Z}}  $ is a stationary sequence of martingale differences with respect to $({\mathcal F}_{j,n} )_{j \in {\mathbb Z}}$ where ${\mathcal F}_{j,n} = {\mathcal G}_{n -j +1}$. 

To end the proof of Remark \ref{commentreversed}, it remains to notice that, as in the proof of Theorem  \ref{ThmartW1cond}, Lemma  \ref{keylemmarev} immediately leads to the conclusion of  Theorem  \ref{ThmartW1cond} for the reversed martingale $M_n$  with  ${\mathcal G}_{n+1} $ replacing $\F_0$. $\diamond$

\subsection{Proof of the results of Section \ref{section-cocycle}}  \label{sectproofsection-cocycle}

\subsubsection{Proof of Proposition \ref{unique}.}  The existence of  a $\mu$-invariant probability 
on $\B(X)$ is a well-known fact based on the compactness of $X$, see for instance 
Lemma 2.16 of \cite{BQ-book}. Hence, we shall only prove uniqueness. 

\smallskip

By Lemma \ref{lemma-complete} and the Borel-Cantelli lemma we see that 
for every $x,y\in X$ such that $f(x)=f(y)$, 
\begin{equation}\label{as-contraction}
d(B_{2^n}\cdot x,B_{2^n}\cdot y) \underset{n\to \infty}\longrightarrow 
0\qquad \as
\end{equation} 

Let $\nu$ be a $\mu$-invariant 
probability on $X$. By Lemma 2.17 of \cite{BQ-book}, 
for $\P$-almost every $\omega\in \Omega$, there exists a probability 
$\nu_\omega$ on $X$ such that for every continuous function $\varphi$  on $X$, we have 
\beq\label{lem-BQ}
\int_X\varphi(B_n(\omega)\cdot x)\nu(dx)\underset{n\to \infty}\longrightarrow 
\int_X\varphi(x)\nu_\omega(dx)\, .
\eeq
 Let $\varphi$ be a  continuous function on $X$. 
For every $\f\in F$, pick  $y_\f\in X_\f$. Using \eqref{as-contraction}, combined with Fubini's theorem for $\P\otimes \nu$ and the uniform continuity of $\varphi$, we infer that for $\P$-almost every $\omega\in \Omega$  we have 
$$
| \varphi(B_{2^n}(\omega)\cdot y_\f)- \varphi(B_{2^n}(\omega)\cdot x)|\underset{n\to \infty}\longrightarrow 0 \qquad \mbox{for $\nu$-almost 
every $x\in X_\f$}\, .
$$
Recall that $(X_\f )_{\f\in F}$ forms a partition of $X$ and that, by Lemma \ref{basic-lemma}, $\nu(X_\f)=\frac1{|F|}$ for every $\f \in F$
Hence, by Lebesgue dominated convergence theorem, using \eqref{lem-BQ}, we infer that, for 
$\P$-almost every $\omega\in \Omega$, 
\begin{equation}\label{eq}
\frac1{|F|}\sum_{\f\in F} \varphi(B_{2^n}(\omega)\cdot y_\f)\underset{n\to \infty}\longrightarrow 
\int_X\varphi(x)\nu_\omega(dx)\, .
\end{equation}
Let $\omega$ be such that \eqref{eq} holds. Since $X$ is compact, there exist $(n_k)_{k\ge 1}$ and $(Z_{\f}(\omega))_{\f\in F}$ 
such that $B_{2^{n_k}}(\omega)\cdot y_\f\underset{k\to \infty}\longrightarrow 
Z_\f(\omega)$.

\smallskip

From the previous considerations, we conclude that  for every continuous $\varphi$ on $X$,
$$
\int_X\varphi(x)\nu_\omega(dx) =\frac1{|F|} \sum_{\f\in F}\varphi(Z_\f(\omega))\, .
$$ 
Then, we infer that $\nu_\omega=\frac1{|F|} \sum_{\f\in F}\delta_{Z_\f(\omega)} $, 
hence does not depend on $\nu$. 

\smallskip

The uniqueness of $\nu$ then follows from the fact that for every 
continuous  $\varphi$ on $X$, 
$$ \int_X\varphi(x)\nu(dx)= \int_\Omega \Big( \int_X\varphi(x)\nu_\omega(dx)\Big)d\P(\omega)\, ,
$$ see Lemma 2.19 of 
\cite{BQ-book}. \hfill $\square$ 

\medskip

\subsubsection{\bf Proof of Lemma \ref{basic-lemma}.} Let $\nu$ be a 
$\mu$-invariant probability on $\B(X)$.  Let $\psi\, :\, F\to \R$. We have 
\begin{gather*}
\int_{G\times X}\psi( f(g\cdot x)) \mu(dg)\nu(dx) =\int_X \psi(f(x))\nu(dx)\,  .
\end{gather*}
Hence 
$$
\int_G \left(\sum_{\f\in F} \psi(s(g)\f) \nu(X_\f) \right)\, \mu(dg)= 
\sum_{\f\in F} \psi(\f) \nu(X_\f)\, .
$$
Hence, $\sum_{\f\in F} \nu(X_\f) \delta_\f$ is an $s(\mu)$-invariant 
probability measure on $F$. Since $s(\mu)$ is adapted it follows from 
the Choquet-Deny theorem for compact groups, see for instance \cite{IK},  (notice that for finite groups the proof is quite direct) that the latter measure is the Haar measure on $F$. \hfill $\square$ 

\medskip

\subsubsection{Proof of Lemma \ref{lemma-complete}.}
Of course, it is enough to prove \eqref{complete-int} with $A_{2^n}$ in place of 
$B_{2^n}$. Let $n_0$ be the integer appearing in Definition \ref{weakly-contracting}.

Set $\eta_n:= \sup_{ x\neq y,f(x)=f(y)} \int_G \log \Big( \frac{d(g\cdot x,g\cdot y)}{d(x,y)}
\Big)\, \mu^{*n}(dg)$. Then for every integers $k\ge 1$ and $\ell\in [0,n_0-1]$, 
$\eta_{kn_0+\ell}\le k\eta_{n_0}+\eta_\ell\le k\eta_{n_0}+\max_{0\le m\le n_0}\eta_m$. Since $\eta_{n_0}<0$, 
it follows that $\eta_n<0$ for every $n$ large enough. Hence we may assume that 
$n_0=2^{r_0}$.

\smallskip
Using that $A_0=e$, we have 
\begin{gather*}
 \log\, \Big(\frac{d(A_{2^{n+r_0}}\cdot x,A_{2^{n+r_0}}\cdot 
y)}{d(x,y)}\Big)=\sum_{k=1}^{2^n}
  \log\, \Big(\frac{d(A_{2^{r_0}k}\cdot x,A_{2^{r_0}k}\cdot 
y)}{d(A_{2^{r_0}(k-1)} \cdot x,A_{2^{r_0}(k-1)}\cdot y)}\Big)
\end{gather*}
As mentionned in Remark \ref{comparewithBQ} the right-hand side of 
\eqref{contracting} may be $-\infty$. To take care of that technical matter let 
us introduce some $R>0$ (that clearly exists) such that 
$$
\sup_{ x\neq y,f(x)=f(y)} \int_G\max\left(-R, \log \Big( \frac{d(g\cdot x,g\cdot y)}{d(x,y)}
\Big)\right)\, \mu^{*2^{r_0}}(dg) <0\, .
$$
Let $\mu_0:=\mu^{*2^{r_0}}$. For every $x,y\in X$ and $g\in G$, set 
\begin{equation} \label{defofphiR}
\varphi_R(g,x,y):= \max\left(-R, \log \Big( \frac{d(g\cdot x,g\cdot y)}{d(x,y)}
\Big)\right)\, .
\end{equation}
Set also $V_k:=Y_{2^{r_0}k}\cdots Y_{2^{r_0}(k-1)+1}$.

\smallskip
Then, we have 
\begin{gather*}
\log\, \Big(\frac{d(A_{2^{n+r_0}}\cdot x,A_{2^{n+r_0}}\cdot 
y)}{d(x,y)}\Big)  \le \sum_{k=1}^{2^n} \varphi_R(V_k,A_{2^{r_0}(k-1)}\cdot x,
A_{2^{r_0}(k-1)}
\cdot y)\\
  = \sum_{k=1}^{2^n} \Big(  \varphi_R(V_k,A_{2^{r_0}(k-1)}\cdot x,A_{2^{r_0}(k-1)}
\cdot y)-\int_G \varphi_R\big(g,A_{2^{r_0}(k-1)}\cdot x,A_{2^{r_0}(k-1)}\cdot y\big)\, 
\mu_0(dg) \Big) \\+  \sum_{k=1}^{2^n} 
\int_G \varphi_R\big(g,A_{2^{r_0}(k-1)}\cdot x,A_{2^{r_0}(k-1)}\cdot y\big)\, \mu_0(dg) \\
:= \sum_{k=1}^{2^n} D_{R,k}(x,y)+\sum_{k=1}^{2^n} E_{R,k}(x,y)\\
:= M_{R,2^n}(x,y)+Q_{R,2^n}(x,y)\, .
\end{gather*}

Notice that $(M_{R,n}(x,y))_{n\in \N}$ is a $(\F_{2^{r_0}n})_{n\in \N}$-martingale. 
In order to prove \eqref{complete-int} it is enough to prove that for every $\varepsilon>0$, 
\begin{equation}\label{mart-part}
\sum_{n\ge 0} \sup_{ x\neq y,f(x)=f(y)}\P(|M_{R,2^n}(x,y)|\ge 2^n\varepsilon )<\infty \, ,
\end{equation}  and that there exists $\tilde \ell >0$ such that
\begin{equation}\label{bounded-part}
\sum_{n\ge 0} \sup_{ x\neq y,f(x)=f(y)}\P\, \left( Q_{R,2^n}(x,y)\ge -2^n\tilde \ell\right)\, <\infty\, .
\end{equation}

We first prove \eqref{bounded-part}. Set 
$$-2\tilde \ell:= \sup_{ x\neq y,f(x)=f(y)} \int_G\max\left(-R, \log \Big( \frac{d(g\cdot x,g\cdot y)}{d(x,y)}
\Big)\right)\, \mu^{*2^{r_0}}(dg) <0\, .$$

Since $Q_{R,2^n}(x,y)\le -2^{n+1}\tilde \ell$, $\P(\, \left( Q_{R,2^n}(x,y)\ge -2^n\tilde \ell\right)=0$ and the result follows.

\smallskip

 Let us prove \eqref{mart-part} (for every $\varepsilon>0$). Define a function 
 $\Lambda$ on $G$ by setting  $\Lambda(g):= \max(R,\log^+{\rm Lip}
(g))$ for every $g\in G$. 

\smallskip

Let $\Gamma_n:=\cup_{k=1}^{2^n} 
\{\Lambda(V_k)> \varepsilon 2^n\}$ and 
\begin{gather*}
D'_{R,k}(x,y):=
D_{R,k}(x,y) {\bf 1}_{\{\Lambda(V_k)\le \varepsilon 2^n\}}\\
D''_{R,k}(x,y) :=D'_{R,k}(x,y)-\E(D'_{R,k}(x,y) |\F_{2^{r_0}(k-1)})\, .
\end{gather*}

Notice that $|D_{R,k}|\le \Lambda(V_k)+\E(\Lambda(V_0))$ and since $D_{R,k}$ is such that $\E(D_{R,k}(x,y)|\F_{2^{r_0}(k-1)}) =0$, 
\[
\E(D'_{R,k}(x,y)|\F_{2^{r_0}(k-1)}) = - \E(D_{R,k}(x,y) {\bf 1}_{\{\Lambda(V_k) > \varepsilon 2^n\}} |\F_{2^{r_0}(k-1)}) \, .
\]
Hence, using independence, 
\begin{multline*}
|\E(D'_{R,k}(x,y)|\F_{2^{r_0}(k-1)})|  \leq \E( \Lambda(V_k) {\bf 1}_{\{\Lambda(V_k) > \varepsilon 2^n\}} |\F_{2^{r_0}(k-1)})  + \E(\Lambda(V_0))  \E(  {\bf 1}_{\{\Lambda(V_k) > \varepsilon 2^n\}} |\F_{2^{r_0}(k-1)})  \\
= \E( \Lambda(V_0) {\bf 1}_{\{\Lambda(V_0) > \varepsilon 2^n\}} )   + \E(\Lambda(V_0)) { \mathbb P}( \Lambda(V_0) > \varepsilon 2^n )  = 3 \E( \Lambda(V_0) {\bf 1}_{\{\Lambda(V_0) > \varepsilon 2^n\}} ) \, .
\end{multline*}
Since $\Lambda\in L^1(\mu_0)$, it follows
\begin{equation}
\label{cond-bound}
|\E(D'_{R,k}(x,y)|\F_{2^{r_0}(k-1)})|\le 3 \int_G \Lambda(g)
{\bf 1}_{\{\Lambda(g)> \varepsilon 2^n\}}
\, \mu_0(dg)\underset{n\to \infty}\longrightarrow 0\, .
\end{equation}
Now we have
\begin{align*}
\P(|M_{R,2^n}(x,y)|\ge 2^n\varepsilon ) & \le \P(\Gamma_n) + 
\P(|\sum_{k=1}^{2^n} D'_{R,k}(x,y)|\ge 2^n\varepsilon ) \\
&  \le 2^n\P(\Lambda(V_0) \ge 2^n\varepsilon )+ 
\P \big ( \big |\sum_{k=1}^{2^n} D''_{R,k}(x,y) \big  |\ge 2^{n-1}\varepsilon \big )\\
& \qquad \qquad  +\P \big  (\big  |\sum_{k=1}^{2^n} \E(D'_{R,k}(x,y)|\F_{2^{r_0}(k-1)}) \big  | \ge 2^{n-1}\varepsilon \big  ) \, .
\end{align*}
Using \eqref{cond-bound}, we see that for $n$ large enough, 
$\P(|\sum_{k=1}^{2^n} \E(D'_{R,k}(x,y)|\F_{2^{r_0}(k-1)}) | \ge 2^{n-1}\varepsilon )=0$.
Hence, for $n$ large enough
\begin{gather*}
\P(|M_{R,2^n}(x,y)|\ge 2^n\varepsilon ) \le  2^n\P(\Lambda(V_0) \ge 2^n\varepsilon )+ \frac1{\varepsilon ^2 2^{n-2}}\E(\Lambda(V_0)^2{\bf 1}_{\{
\Lambda(V_0)\leq  2^n\varepsilon\}})\, .
\end{gather*}
Then, \eqref{mart-part} follows from standard computations and the fact that $\E(\Lambda(V_0)) <\infty$. \hfill $\square$

\medskip

\subsubsection{Proof of Lemma \ref{lem-comp}.}  Let $n_0$ be as in Definition \ref{weakly-contracting}. 
%By the argument used in the above remark, it is enough to prove that 
%there exists $\ell>0$ such that 
%\begin{equation*}
%\sum_{k\ge 1} k^{p-2}\sup_{ x, y\in X, x\neq  y} \P \left(
%\max_{1\le j\le kn_0} \log \left(d(A_{j-1}\cdot  x, 
%A_{j-1}\cdot  y)\right)\ge -\ell k\right)<\infty\, .
%\end{equation*}

Using that $(Y_n)_{n\ge 1}$ is iid we see that for every $k\le j\le 2k$ 
$$\sup_{ x, y\in X, x\neq  y} \P \left(
 \log \left(d(A_{j-1}\cdot  x, 
A_{j-1}\cdot  y)\right)\ge -\ell k\right) \le \sup_{ x, y\in X, x\neq  y} \P \left(
 \log \left(d(A_{k-1}\cdot  x, 
A_{k-1}\cdot  y)\right)\ge -\ell k\right)\, .
$$
Similarly, for every $mn_0\le k\le (m+1)n_0$, for some $m\ge 3$, we have (notice that $mn_0-1\ge (m-1)n_0$ and $m+1\le 2(m-1)$)
\begin{gather*}
\sup_{ x, y\in X, x\neq  y} \P \left(
 \log \left(d(A_{k-1}\cdot  x, 
A_{k-1}\cdot  y)\right)\ge -\ell k\right)\\ \le \sup_{ x, y\in X, x\neq  y} 
\P \left( \log \left(d(A_{(m-1)n_0}\cdot  x, 
A_{(m-1)n_0}\cdot  y)\right)\ge - 2 \ell (m-1)n_0\right)\\
\end{gather*}
Let $k\ge 1$. Proceeding as in the proof of Lemma \ref{lemma-complete}, we see that, 
setting $W_k:=Y_{kn_0}\cdots Y_{(k-1)n_0+1}$ and using the notation \eqref{defofphiR}, 
\begin{gather*}
\log \left(d(A_{kn_0}\cdot  x, A_{kn_0}\cdot  y)\right)- \log \left(d(  x,   y)\right)
\\ \le\sum_{j=1}^{k}\Big(\varphi_R(W_j,A_{(j-1)n_0}
\cdot x,A_{(j-1)n_0}\cdot y) -\int_G \varphi_R(g,A_{(j-1)n_0}
\cdot x,A_{(j-1)n_0}\cdot y)\mu^{*n_0}(dg)\\\qquad \qquad \qquad  + 
\sum_{j=1}^k \int_G \varphi_R(g,A_{(j-1)n_0}
\cdot x,A_{(j-1)n_0}\cdot y)\mu^{*n_0}(dg)\\
:=M_k+R_k\, .
\end{gather*}
Clearly, $R_k\le -\gamma k$ where $\gamma:=-\sup_{x,y\in X,x\neq y}\int_G \varphi_R(g,A_{(j-1)n_0}
\cdot x,A_{(j-1)n_0}\cdot y)\mu^{*n_0}(dg)>0$. Hence, 
\eqref{complete-cocycle} will hold, say with $\ell=\gamma/(4n_0)$, if we can prove  prove that  for every $\varepsilon>0$,
\begin{equation}\label{comp-aux}
\sum_{k\ge 1} k^{p-2}\P \Big (\max_{1\le j\le k}|M_j|> \varepsilon k \Big )<\infty\, .
\end{equation}
%Notice that \eqref{comp-aux} without the maximum inside the probability would be sufficient to get  \eqref{complete-cocycle}, but we shall need the maximum to derive \eqref{comp-aux-2} below. 

Now, the fact that \eqref{comp-aux} holds for every $\varepsilon>0$  may be proved as Theorem 4.1 of \cite{CDM}, where a more general result is obtained.

\smallskip

Hence, the proof of \eqref{complete-cocycle} is completed and it remains to prove \eqref{exp-conv}.

\smallskip 

Taking $p=1$ in  \eqref{comp-aux} and using that, for $2^\ell\le k\le 2^{\ell+1}-1$, $\P(\max_{1\le j\le 2^\ell}|M_j|> \varepsilon 2^{\ell})
\le \P(\max_{1\le j\le k}|M_j|> \varepsilon k/2)$, we infer that 
\begin{equation*}\label{comp-aux-2}
\sum_{k\ge 1} \P \Big (\max_{1\le j\le 2^k}|M_j|> \varepsilon 2^{k} \Big )<\infty\, ,
\end{equation*}
for every $\varepsilon>0$ and \eqref{exp-conv} follows from the Borel-Cantelli lemma, since 
$R_k\le -\gamma k$.\hfill $\square$

%Let us prove that there exists $\tilde \ell$, such that \eqref{bounded-part} 
%holds. We proceed as in \cite{CDJ} (proof of Lemma 6). The main idea 
%of the proof seems to be due (at least in this context) to Lepage \cite{lepage}.

%\smallskip

%Notice that 
%$$
%|E_{R,k}(x,y)|\le K:= \int_G\Delta (g)\mu(dg)\, .
%$$
%Hence, for every $n\in \N$, and every $a>0$ we have 
%$$
%|\E( {\rm e}^{aE_{R,k}(x,y)}-1-a \E
%$$

\subsubsection{Proof of Theorems \ref{ThASIPcocycles} and \ref{ThW1condcocycles}.} 

Let $n\ge 1$. We have 
\begin{equation}\label{Tn}
T_n:=\sigma(A_n,W_0)-n \lambda_\mu=\sum_{k=1}^n (\sigma(Y_k,A_{k-1}\cdot W_0)-\lambda_\mu)\, ,
\end{equation}
where all the summands have the same law. 

\smallskip

We start by noticing that, by Proposition \ref{prop-aperiodic}  (and Remark \ref{remark-cocycle}), using that $p\ge 2$, 
\begin{align}\label{serie-cobord}
\sum_{n\ge 1} \|\E\big(\sigma(Y_n,A_{n-1}\cdot W_0\, |\sigma\{W_0\}\big)-\lambda_\mu\|_p \le 
\sum_{n\ge 1} n^{p-2}\|\E\big(\sigma(Y_n,A_{n-1}\cdot W_0\, |\sigma\{W_0\}\big)-\lambda_\mu\|_p  \nonumber  \\
\le \sum_{n\ge 1} n^{p-2} \sup_{x \in X} \vert \E\big(\sigma(Y_n,A_{n-1} x \big)-\lambda_\mu \vert <\infty \, .
\end{align}
 Hence, by Gordin's $L^p$-criteria, we infer  that 
$\sigma(Y_n,A_{n-1}\cdot W_0)-\lambda_\mu= D_n+ R_n-R_{n-1}$, where $(D_n)_{n\ge 1}$ is a stationary sequence of martingale differences in $L^p$ and $(R_n)$ is a stationary process in $L^p$. Write $M_n:= \sum_{k=1}^n D_k$ and $(M_{n})_i:= \langle M_n, e_i\rangle$ for $1\le i\le d$.

\smallskip

Then, it is clearly sufficient to prove the theorems for $(M_n)_{n\ge 1}$ instead of $(T_n)_{n\ge 1}$. Since we are back to the study of a 
martingale with stationary increments, we only have to check the conditions of Theorems  \ref{Thmart} and \ref{ThmartW1cond}. In particular, it suffices to check that
%, using the notations of Remark \ref{remark-cocycle} 
\beq \label{first-cond}
\sum_{n=1}^{\infty}\frac{1}{n^{3-p/2}}\big \Vert \bkE \big
((M_n)_i (M_n)_j\big | \F_{0} \big ) - \bkE \big
((M_n)_i (M_n)_j \big ) \big \Vert_{p/2} < \infty  \, .
\eeq

Clearly, using the basic equality $(a+b)^2 - a^2 - b^2 = 2 ab$, the latter holds provided that one can prove that 
\[ 
\sum_{n=1}^{\infty}\frac{1}{n^{3-p/2}}\big \Vert \bkE \big
({\widetilde M}_n^2\big | \F_{0} \big ) - \bkE \big
({\widetilde M}_n^2 \big ) \big \Vert_{p/2} < \infty \, ,
\]
for any ${\widetilde M}_n\in \{(M_n)_i\, :\, 1\le i\le d\}\cup\{(M_n)_i+(M_n)_j\, :\, 1\le i<j\le d\}$. In particular, we are back to prove that \eqref{first-cond} holds in the case of an $\R$-valued cocycle. As in the proof of Propositions 8 and 9 of \cite{CDJ}, taking into account 
\eqref{serie-cobord}, it is enough to prove that 
\beq \label{cond-Tn1}
\sum_{n=1}^{\infty}\frac{1}{n^{3-p/2}}\big \Vert \bkE \big
(T_n^2\big | \F_{0} \big ) - \bkE \big
(T_n^2 \big ) \big \Vert_{p/2} < \infty \, .
\eeq

Then, proceeding as in the proof of Proposition 12 of \cite{CDJ} (based on Proposition 11 there), making use of \eqref{simple-series}, \eqref{square-series} 
and \eqref{mixed-series}, we infer that \eqref{cond-Tn1}   holds. \hfill $\square$

%\begin{Theorem}
%Let $p\ge 1$. Assume that the action is weakly $(\mu,p)$-contracting, that 
%$\mu$ is $F$-strictly aperiodic and that  
%$\sigma$ admits  a moment of order $p$. We have 
%\begin{itemize}
%\item [$(i)$] Items $(i)$ to $(iv)$ in  Theorem 1 of \cite{CDJ}, as well as 
%items $(i)$ and $(ii)$ in Theorem 2 of \cite{CDJ}, hold with the corresponding $p$ 
%if one replaces $S_{n,\overline x}$ with $\sigma(A_n,x)$;
%\item [$(ii)$] Theorem 4.1 of "deviations" holds for the corresponding $p$ and, 
%if $p=2$, Theorem 4.2 holds, if one replaces  $\log\|A_n x\|$ with $\sigma(A_n,x)$.
%\end{itemize}
%\end{Theorem}
%
%\begin{Theorem}
%Assume that the action is weakly $(\mu,2)$-contracting, that 
%$\mu$ is $F$-strictly aperiodic  and that there exists 
%$r,\delta >0$ such that 
%$$
%\int_G{\rm e}^{\delta \kappa_0^r(g)}\mu(dg) <\infty\, .
%$$  
%Then, if $r\ge 1$ (resp. $r\in (0,1)$), Theorem 2.1 (resp. Theorem 2.2) of "deviations" holds if one replaces  $\log\|A_n x\|$ with $\sigma(A_n,x)$.
%\end{Theorem}
%
%{\bf Mentionner egalement le theoreme 2.3}. 

\subsubsection{Proof of Theorems \ref{ThASIP-Iwasawa} and \ref{ThW1-Iwasawa}.} The case where $V_n=\sigma(A_n,W_0)$ follows directly from Theorems \ref{ThASIPcocycles} and \ref{ThW1condcocycles}.  It remains to prove them for the other possible values of $V_n$. We shall only complete the proof of Theorem \ref{ThASIP-Iwasawa} since the proof of 
Theorem \ref{ThW1-Iwasawa} may be done as the proof of item $2.$ of Theorem  \ref{ThASIP-Iwasawa}. The proofs make use of results  from \cite{BQ-book}.

\smallskip

Let us  start by proving Item $1$ of Theorem \ref{ThASIP-Iwasawa}.  It follows from item $(d)$ of Theorem 10.9 of \cite{BQ-book}  that 
$\P$-almost surely
$$
\sup_{n\ge 1} |\sigma(A_n, W_0)-\kappa(A_n)|_\mfa<\infty\, .
$$
Hence Item $1$ for the Cartan projection follows from the result for $(\sigma(A_n,W_0) - n \lambda_{\mu} )_{n \geq 1}$. Similary, Theorem 10.9 of \cite{BQ-book} implies that for every $x\in {\mathcal P}$, 
$$
\sup_{n\ge 1} |\sigma(A_n,x)-\kappa(A_n)|_\mfa < \infty\, .
$$
Hence Item $1$ for $(\sigma(A_n, x) - n \lambda_{\mu})_{n\ge 1}$, for every $x\in {\mathcal P}$, follows from the result for $(\kappa(A_n) - n \lambda_{\mu} )_{n \geq 1}$.
\smallskip

Let us prove Item $2$ of Theorem \ref{ThASIP-Iwasawa}. By (13.31) of \cite{BQ-book}, there exist $M>0$ and $\eta_1, \ldots , \eta_r\in {\mathcal P}$, 
such that, for every $g\in G$ there exists $i\in \{1, \ldots, r\}$ such that
\begin{equation}\label{sigma-kappa}
|\sigma(g, \eta_i) -\kappa(g)|_\mfa \le M\, .
\end{equation}
Let $(\Gamma_i)_{1\le i\le r}$ be a partition of $\Gamma_\mu$ such that for every $g\in \Gamma_i$ \eqref{sigma-kappa} holds. Let $(N_n)_{n\ge 1}$ be the sequence of iid normal variables appearing in Item $2$ and approximating $(\sigma(A_n, W_0))_{n\ge 1}$. Set $W_n:=\sum_{k=1}^n N_k$. We have 
\begin{gather*}
\big| \kappa(A_n)-W_n \big|_\mfa=\sum_{\f\in F} \quad \sum_{i\in \{1,\ldots ,r\},:\, \eta_i\in X_\f} \big| \kappa(A_n)-W_n\big|_\mfa 
 {\bf 1}_{\Gamma_i}(A_n)\\
\le M+  \sum_{\f\in F} \quad \sum_{i\in  \{1,\ldots ,r\},:\, \eta_i\in X_\f}\big| \sigma(A_n,\eta_i)- W_n|_\mfa {\bf 1}_{\Gamma_i}(A_n)\, .
\end{gather*}
For $\eta_i\in X_\f$ write
$$
\sigma(A_n,\eta_i)- W_n= \frac1{|F|}\Big[ \int_{X_\f} \big( \sigma(A_n,\eta_i)-\sigma(A_n,x)\big)\nu(dx) + 
\int_{X_\f} \big( \sigma(A_n,x)-W_n\big)\nu(dx)\Big]\, .
$$ 
Combining the above estimates and using Proposition \ref{startingpoint}, we infer that there exists $K>0$ 
$$
\big \|\, \big| \kappa(A_n)-W_n \big|_\mfa\, \big\|_1\le rK+ r\|\, |\sigma(A_n,W_0)-W_n|_\mfa\, \|_1 \, ,
$$
and the desired result follows.

\medskip

The case where $G$ is connected follows from Lemma \ref{startingpoint} (see the remark after the lemma). The case where $G=GL_d(\R)$ follows also from 
Lemma \ref{startingpoint} once one has noticed that if 

$$\varepsilon=\left(
\begin{array}{ccccc}
-1  &  0  &  \ldots &  \ldots  &  0\\
0  &  1  &  0  &  \ldots  &  0\\
\vdots &\ddots   &  \ddots &  \ddots  &  \vdots\\
\vdots  &  & \ddots  &  \ddots  &  0 \\
0 &  \ldots & \ldots  & 0  &  1
\end{array}
 \right)
$$
then $f(kP_c)=-f(k\varepsilon P_c)$ and $\sigma(g, kP_c)=\sigma(g, k\varepsilon P_c)$ for every $g\in GL_d(\R)$ and every $kP_c\in G/P_c$. \hfill $\square$

\section{Appendix}

\setcounter{equation}{0}

The following result is a Fuk-Nagaev type inequality for martingales with moments of order $p \in ]2,3]$. 

\begin{Proposition} \label{inegamax}
Let $(d_i)_{i \in \mathbb Z}$ be a real stationary martingale differences
sequence with respect to $({\cal F}_i)_{i \in {\mathbb Z}}$. Let
$\sigma$ denote the standard deviation of $d_0$. Let $p \in ]2,3]$.
Assume that $\bkE |d_0 |^p < \infty$ and that the conditions
(\ref{C1}) and (\ref{C2}) are fulfilled. Let $M_n^* = \max_{1 \leq k
\leq n} M_k$. Then for any positive real $x$,
$$
\BBP (M_n^* \geq x) \leq 7 \times 2^{2p-1/2} \Big (
\frac{n\sigma^2}{x^{2}} \Big )^{p+1/2}\exp \Big( - \frac{x^2}{8n
\sigma^2} \Big ) + C_p n x^{-p} \, ,
$$
where $C_p$ is a positive constant depending only on $p$,  $\bkE |d_0 |^p$ and on the two series (\ref{C1}) and (\ref{C2}) but not on $n$. 
\end{Proposition}

\noindent{\bf Proof of Proposition \ref{inegamax}.} The proof uses similar arguments as in the proof of Proposition 5.2 in \cite{MR}.

Let $x$ be a
positive real and $\varphi$ be the function from $\mathbb R$ to
${\mathbb R}^+$ defined by
$$
\varphi(t) = \frac{(t -x/2)^{p}_+}{p(p-1)} \, .
$$
Since $\varphi$ is nonnegative and convex, $\varphi(M_n)$ is a
submartingale. Consequently Doob's maximal inequality entails that
$$
\BBP (M_n^* \geq x)  = \BBP \big (\varphi(M_n^*) \geq \varphi(x)
\big ) \leq \frac{\bkE(\varphi(M_n))}{\varphi(x)} \, .
$$
Let now $Y$ be a standard gaussian random variable independent of $(d_i)_{i \in \mathbb Z}$. Then
\beq \label{p1inemax} \BBP (M_n^* \geq x)   \leq \frac{2^p p
(p-1)}{x^p} \big ( \bkE(\varphi(\sigma \sqrt n Y)) + \bkE(\varphi(M_n)) - \bkE(\varphi(\sigma \sqrt n Y))  \big )  \, . \eeq Since $\varphi''(t)=(t
-x/2)^{p-2}_+$, it follows that $\varphi''$ is $(p-2)$ H\"older.
Then,
$$
\bkE(\varphi(M_n)) - \bkE(\varphi(\sigma \sqrt n Y)) \leq \sup_{f
\in {\Lambda_p}}\bkE( f(M_n) - f(\sigma \sqrt n Y)) \, ,
$$
where (recall that $p \in ]2,3]$) $\Lambda_p$ is the class of real
functions $f$ which are $2$-times continuously differentiable and
such that $$ |f^{''}(x) - f^{''}(y) | \leq | x - y |^{p-2} \ \hbox{
for any } (x,y) \in \bkR \times \bkR \, . $$ Now denoting by
$P_{M_n}$ the law of $M_{n}$ and  by $G_{n\sigma^2}$  the normal
distribution $ N(0, n\sigma^2)$ we have
$$
\sup_{f \in {\Lambda_p}}\bkE( f(M_n) - f(\sigma \sqrt n Y)) :=
\zeta_p ( P_{M_n}, G_{n \sigma^2}) \, ,
$$
that is the Zolotarev distance of order $p$ between $P_{M_n}$  and
 $G_{n\sigma^2}$. Applying Theorem 2.1 in \cite{DMR09}, we derive that, under (\ref{C1}) and (\ref{C2}), there
exists a positive constant $K_p$ depending  only on $p$,  $\bkE |d_0 |^p$ and on the two series (\ref{C1}) and (\ref{C2}), but not on $n$, such that \beq \label{majzetap} \bkE(\varphi(M_n)) -
\bkE(\varphi(\sigma \sqrt n Y))  \leq K_p n \, . \eeq On an other
hand, using the fact that $(u+x/2)^2\ge ux+x^2/4$, we
get that
$$\bkE(\varphi(\sigma \sqrt n Y)) \leq \frac{e^{-x^2/(8n\sigma^2)}}{p(p-1)}
\int_{0}^{\infty} e^{-xu/(2n \sigma^2)} \frac{u^p}{\sigma \sqrt{2 n
\pi}} du \, .
$$With the change of variables $v= xu/(2n\sigma^2)$ and
since $p\in ]2,3]$,  we derive that
\begin{eqnarray} \label{majgauss} p(p-1)\bkE(\varphi(\sigma \sqrt n Y))
&\leq& e^{-x^2/(8n\sigma^2)} \Big ( \frac{2n\sigma^2}{x}\Big )^{p+1}
\frac{1}{\sigma \sqrt{2 n \pi}}\int_0^{\infty} v^p e^{-v} dv
\nonumber \\ &\leq & \frac{7}{2}
\frac{(2n\sigma^2)^{p+1/2}}{x^{p+1}} e^{-x^2/(8n\sigma^2)} \, .
\end{eqnarray} Starting from (\ref{p1inemax}) and using
(\ref{majzetap}) and (\ref{majgauss}) we get the result.

\bigskip

\noindent{\bf Acknowledgements} We would like to thank \c{C}a$\breve{\text g}$ri Sert for useful discussions and references on the Iwasawa cocycle.

\end{document}